\documentclass[11pt,reqno]{amsart}

\usepackage[a4paper]{geometry}
\usepackage{float}
\usepackage[usenames,dvipsnames]{xcolor}
\usepackage{subfigure}

\usepackage[colorlinks, linkcolor=blue,anchorcolor=Periwinkle,
   citecolor=red,urlcolor=Emerald]{hyperref}
\usepackage{amsmath}           
\usepackage{amssymb}           
\usepackage{latexsym}          
\usepackage{bookmark}
\usepackage{enumitem}
\usepackage{array}
\usepackage{longtable}
\usepackage{graphicx}
\usepackage[bottom]{footmisc}
\usepackage{cleveref}
\usepackage{transparent}

\usepackage{tikz}
\usetikzlibrary{matrix}
\usepackage{url}

\usetikzlibrary{decorations.markings,cd}
\usetikzlibrary{arrows}
\usetikzlibrary{decorations.pathreplacing,decorations.pathmorphing,shadings,fadings,calc}
\usepackage{mathtools}
\IfFileExists{extarrows.sty}{\usepackage{extarrows}}{}

\providecommand{\xlongleftrightarrow}[2][]{\xleftrightarrow[#1]{#2}}
\makeatletter
\providecommand{\reviewequalfill}{\arrowfill@=\relbar=}
\providecommand{\xlongequal}[2][]{\ext@arrow 0055\reviewequalfill{#1}{#2}}
\makeatother
\IfFileExists{stmaryrd.sty}{\usepackage{stmaryrd}}{}

\tikzset{->-/.style={decoration={  markings,  mark=at position #1 with
    {\arrow{>}}},postaction={decorate}}}
\tikzset{-<-/.style={decoration={  markings,  mark=at position #1 with
    {\arrow{<}}},postaction={decorate}}}
\tikzcdset{arrow style=tikz, diagrams={>=stealth}}

\newcolumntype{L}{>{$}l<{$}} 

\usepackage{tikz}
\usetikzlibrary{matrix, calc, arrows,cd,
decorations.pathreplacing, decorations.markings, decorations.pathmorphing,
backgrounds,fit,positioning,shapes.symbols,chains,shadings,fadings}

\IfFileExists{mathabx.sty}{\usepackage{mathabx}}{%
\newcommand{\Neptune}{\mathord{\text{\tikz[baseline=-.5ex,scale=.13]{\draw[line width=.45pt] (0,0)--(0,1.8) (-.6,1.8)--(-.6,.8)--(.6,.8)--(.6,1.8) (-.35,.35)--(.35,.35);}}}}
}

\theoremstyle{plain}
\newtheorem{theorem}{Theorem}[section]
\newtheorem{lemma}[theorem]{Lemma}
\newtheorem{corollary}[theorem]{Corollary}
\newtheorem{proposition}[theorem]{Proposition}

\newtheorem{cordef}[theorem]{Definition/Corollary}
\newtheorem{thmdef}[theorem]{Definition/Theorem}
\newtheorem{conjecture}[theorem]{Conjecture}
\newtheorem{assumption}[theorem]{Assumption}
\newtheorem{convention}[theorem]{Convention}

\theoremstyle{definition}
\newtheorem{definition}[theorem]{Definition}

\newtheorem{remark}[theorem]{Remark}

\numberwithin{equation}{section}

\newcommand\hua{\mathcal}

\newcommand\ZZ{\mathbb{Z}}
\newcommand\QQ{\mathbb{Q}}
\newcommand\RR{\mathbb{R}}

\renewcommand{\setminus}{\smallsetminus}

\newcommand{\rank}{\operatorname{rank}}
\newcommand{\id}{\operatorname{id}}

\newcommand\Add{\operatorname{add}}
\newcommand\Ind{\operatorname{Ind}}

\newcommand\Hom{\operatorname{Hom}}

\newcommand\Ext{\operatorname{Ext}}

\newcommand{\C}{\operatorname{\hua{C}}} 
\newcommand{\per}{\operatorname{per}} 

\newcommand{\CEG}{\operatorname{CEG}} 
\newcommand{\uCEG}{\underline{\CEG}} 

\renewcommand{\k}{\mathbf{k}}

\newcommand{\on}[1]{\operatorname{#1}}

\usepackage{pifont}

\def\be{\begin{equation}}
\def\ee{\end{equation}
}

\def\CCx{\on{Cpx}}   
\def\cc{\on{cell}}  
\def\dcc{\overrightarrow{\on{cell}}}
\def\sphere{\on{S}}   
\def\susp{\Sigma}       

\def\ww{\mathbf{w}}
\def\vv{\mathbf{v}}
\def\uu{\mathbf{u}}

\newcommand{\DynA}[1]{A_{#1}}

\def\source{X[1]}

\newcommand{\flowd}[1]{\on{evo}^{#1}_{\downarrow}}

\newcommand{\flowu}[1]{\on{evo}^{#1}_{\uparrow}}

\newcommand{\flowdu}[1]{\on{evo}_{\updownarrow}^{#1}}

\newcommand{\Dflow}[1]{\on{Evo}^{\downarrow}_{#1}}
\newcommand{\Uflow}[1]{\on{Evo}^{\uparrow}_{#1}}
\newcommand{\DUflow}[1]{\on{Evo}^{\updownarrow}_{#1}}

\newcommand{\flowD}[1]{\on{Evo}_{\downarrow}^{#1}}
\newcommand{\flowU}[1]{\on{Evo}_{\uparrow}^{#1}}
\newcommand{\flowDU}[1]{\on{Evo}_{\updownarrow}^{#1}}

\newcommand{\pcts}[1]{\mathbf{#1}}
\newcommand{\pcto}[1]{\mathbf{#1}_\oplus} 
\newcommand{\ccpt}[1]{\mathrm{#1}}
\newcommand{\real}[1]{\underline{#1}}

\def\xevo{\ccpt{X}}

\def\con{\sigma}
\def\fan{\Delta}

\newcommand{\bfit}[1]{\textbf{\textit{#1}}}

\newcommand{\ubf}[1]{\underline{\mathbf{#1}}}

\def\PP{\mathbb{P}}

\def\seed{\Sigma}
\def\seeds{\mathbf{\Sigma}}
\def\xx{\mathbf{x}}
\def\yy{\mathbf{y}}
\def\zz{\mathbf{z}}
\newcommand{\pos}[1]{[#1]_+}
\newcommand{\nega}[1]{[#1]_-}
\def\regT{\mathbb{T}}

\def\cvec{\bfit{c}}
\def\gvec{\bfit{g}}
\def\Cmat{{C}}
\def\Gmat{{G}}
\def\Fpoly{{F}}

\def\coex{a}

\def\supp{\on{supp}}

\newcommand{\Bgzp}[2]{\operatorname{Bgz}^{\uparrow}_{#2}(#1)}
\newcommand{\Bgzm}[2]{\operatorname{Bgz}^{\downarrow}_{#2}(#1)}
\newcommand{\Bgzpm}[2]{\operatorname{Bgz}^{\updownarrow}_{#2}(#1)}

\newcommand{\tBgzpm}[2]{\widetilde{\operatorname{Bgz}}^{\updownarrow}_{#2}(#1)}

\newcommand{\suppd}[1]{\operatorname{\Neptune}_{\downarrow}^{#1}}
\newcommand{\suppu}[1]{\operatorname{\Neptune}_{\uparrow}^{#1}}

\newcommand{\suppdu}[1]{\operatorname{\Neptune}_{\updownarrow}^{#1}}
\newcommand{\ol}[1]{\overline{#1}}

\newcommand{\DT}[1]{#1^{\vee}}
\newcommand{\iDT}[1]{#1^{-\vee}}
\def\Xevo{\DT{\ccpt{X}}}

\def\sCpx{\RR\CCx}
\def\dCpx{\overrightarrow{\CCx}}
\def\ZZp{\ZZ_{\ge0}}
\def\ZZm{\ZZ_{\le0}}

\newcommand{\dom}[1]{\mathcal{T}^{#1}}
\newcommand{\rdom}[1]{\widetilde{\mathcal{T}}^{#1}}

\def\normal{\on{N}}
\def\csum{\epsilon}
\def\tgvec{\widetilde{\gvec}}

\title[$\mathrm{X}$-evolution flows/foliations: supports and orders]
{$\mathrm{X}$-evolution flows/foliations for cluster algebras:\\
partial Bongartz (co-)completions and acyclicity}

\author{Yu Qiu}
\address{Qy: Yau Mathematical Sciences Center and Department of Mathematical Sciences, Tsinghua University, 100084 Beijing, China.
    \&
Beijing Institute of Mathematical Sciences and Applications, Yanqi Lake, Beijing, China}
\email{yu.qiu@bath.edu}

\author{Liheng Tang}
\address{Tl: Qiuzhen College, Tsinghua University, 100084 Beijing, China.}
\email{tang-lh20@mails.tsinghua.edu.cn}

\thanks{Qy is supported by National Natural Science Foundation of China (Grant No. 12425104).}
\dedicatory{Dedicated to Bernhard Keller on the occasion of his sixty-fifth birthday}

\newcommand{\reviewfigure}[2][]{%
\IfFileExists{#2}{\includegraphics[#1]{#2}}{%
\fbox{\parbox[c][2cm][c]{.75\linewidth}{\centering
Figure file not supplied: \texttt{\detokenize{#2}}}}}}

\begin{document}
\def\EVO{evo }
\begin{abstract}

We first generalize the notion of $\xevo$-foliation on cluster complexes to the context of cluster algebras.
This is a continuous generalization of Keller's green mutation, which was previously introduced in the prequel within a categorical setting.
We show that the wall-crossing structure (of possibly higher co-dimension) is encoded by $\xevo$-supports, which form a partial version of Cao-Gyoda-Yurikusa's Bongartz (co-)completions for cluster algebras.
Moreover, we prove that for any cell transversed nontrivially by a fixed $\xevo$-leaf, its rank is an invariant bounded below by the rank of $\xevo$.
Finally, we prove the acyclicity conjecture for $\xevo$-foliation, which also implies the compactness conjecture for any $\xevo$-foliation in finite type cases.


\end{abstract}
\maketitle

\tableofcontents

\addtocontents{toc}{\protect\setcounter{tocdepth}{1}}
\setlength\parindent{0pt}
\setlength{\parskip}{5pt}
\setlength{\emergencystretch}{2em}

\section{Introduction}
\subsection{Cluster algebras and green mutation}\
Cluster algebras, introduced by Fomin and Zelevinsky \cite{FZ,FZ03,FZ07}, have become a key theory because they bridge many areas of mathematics. For example,
\begin{itemize}
    \item Lie theory and total positivity \cite{BFZ},
    \item representation theory and categorification \cite{BMRRT,Ke4},
    \item surface topology and Teichm\"uller theory \cite{FST,FG},
    \item DT-theory \cite{KS,Nag,Ke2}, and
    \item mirror symmetry and tropical geometry \cite{GHKK18}.
\end{itemize}
Over the last two decades, many original conjectures have been proved, often with unexpected tools and techniques. Hence the theory is still active and thriving.

Quiver mutation is a key ingredient in cluster theory.
The cluster exchange graph $\uCEG(\seeds)$ of a cluster algebra/pattern $\seeds$ has reachable clusters as vertices and mutations as edges. 
Since mutation is an involution, $\uCEG(\seeds)$ is unoriented. 
The cluster complex $\CCx(\seeds)$ of $\seeds$ is a simplicial complex whose simplices are partial clusters. 
Thus $\uCEG(\seeds)$ is the dual graph of $\CCx(\seeds)$ with respect to its maximal cells and their codimension-one faces.
Cluster exchange graphs/complexes are generalizations of associahedra and complexes of marked surfaces, cf. \cite{FZ03,FST}. 
They play central roles in applying cluster theory. 
Namely, the former is closely related to exchange graphs of t-structures/hearts (cf. \cite{KQ1}) and the latter to scattering diagrams (cf. \cite{Nak}).

Once an initial cluster $\xx$ is chosen, a mutation is oriented and becomes the green mutation, in the sense of Keller \cite{Ke2}. 
In other words, there is an induced $\xx$-orientation on $\uCEG(\seeds)$ with a unique source and a unique sink in good case. 
The motivation comes from the study of quantum dilogarithm identities and Donaldson--Thomas theory, cf. \cite{KS,Nag}, cf. \cite{Ke2,Ke3}. 
Green mutation has been heavily studied with many applications, e.g. \cite{BDP} and more recently \cite{Yak}.
In particular, green mutation naturally arises from the partial order on hearts/t-structures, cf. \cite{KQ1,Q1}, with applications to fundamental groups of moduli spaces \cite{KQ2}.

\subsection{$\xevo$-evo flow and $\xevo$-foliation}\

In the prequel \cite{QT}, 
we introduce the \emph{$\xevo$-evolution flow} as a continuous generalization of Keller's green mutation.
For any point $\xevo$ in the cluster complex $\CCx(\C)$, we introduced a direction field, called the $\xevo$-evolution flow, 
with only two singularities.
It turns out that such flows give a continuous generalization of green mutations, cf. \cite[Thm.~5.10]{QT}.
The setting there is categorical and the method is representation theoretical,
cf. \Cref{app:cats} for more details.

In this paper, our first aim is to generalize the $\xevo$-evo flow/foliation to the setting of cluster algebras. 
Inspired by \cite[Thm.~5.10]{QT}, we first introduce the \emph{parallel part of the $\xevo$-evo flow}. For a \emph{generic} (that requires all $a_i>0$) point 
\begin{gather}\label{eq:intro-generic}
    \xevo=\sum_{i=1}^n a_i x_i,\qquad \sum_{i=1}^n a_i=1,
\end{gather}
supported on a cluster $\xx=\supp\xevo=(x_1,\ldots,x_n)$, cf. \eqref{eq:supp}, and set the parallel part of the $\xevo$-evo flow in a cluster $\vv=(v_1,\ldots,v_n)$ as follows.
\begin{equation}\label{eq:intro-parallel}
\begin{aligned}
\flowdu{\xevo}(\vv)
 =\sum_{i=1}^n a_i\left(
   \normal_\pm(\tgvec_{x_i}^{\vv})
  -\normal_\mp(\tgvec_{x_i}^{\vv})\right).
\end{aligned}
\end{equation}
Here:
\begin{itemize}
    \item $\updownarrow=\uparrow/\downarrow$, $\pm=+/-$ and $\mp=-/+$.
    \item
    $\tgvec_{x_i}^{\vv}=(\widetilde{g}_{ij})_{v_j\in\vv}$ is the $x_i$-row of the $C$-matrix $C^\xx_\vv$, cf. \eqref{eq:dual3}.
    \item 
    $\normal_\pm(\cvec)$ is the normalized positive/negative part of a vector $\cvec=(c_{v_i})_{v_i\in\vv}$ supports on $\vv$
    \begin{equation}\label{eq:intro-normal}
            \normal_\pm(\cvec)=
            \begin{cases}
                \displaystyle
                \frac{\sum_j[ c_{v_j} ]_\pm v_j}
                     {\sum_j[ c_{v_j} ]_\pm},
                &\sum_j[ c_{v_j} ]_\pm\ne0,\\[1.2em]
                \bfit{0},
                &\sum_j[ c_{v_j} ]_\pm=0,
            \end{cases}
    \end{equation}
    where $[c]_\pm=\max\{\pm c,0\}$.
\end{itemize}
The point is that 
the signs of the nonzero $v$-coordinates in $\normal_\pm(\tgvec^\vv_{x_i})$ are the signs of the $c$-vectors $\cvec^{\xx}_{v,\vv}$ of $v$.
The direction of the downward (resp. upward) flow is controlled by its parallel part and therefore agrees with the green (resp. red) mutation.

Next we extend the parallel part to the case when $\vv$ is a partial cluster with the following key property:
\begin{itemize}
    \item the flow originates from $\vv$, so it supports on a cluster containing $\vv$;
    \item the direction of downward (resp. upward) flow should agree with the green (resp. red) mutation, so the $c$-vectors of the target cluster variables with respect to $\xx$ is positive (resp. negative).
\end{itemize}
The two conditions together are precisely the definition of Bongartz completions/co-completions in the cluster-algebraic setting introduced by Cao--Gyoda--Yurikusa \cite{CGY}, see \Cref{def:bong}. 
The existence of Bongartz co-completion for all partial clusters is guaranteed by the following mild assumption, which will always be assumed throughout the paper for simplicity.
\begin{assumption}\label{ass0}
    $\seeds$ {admits a green-to-red sequence}. 
\end{assumption}
Let $\Bgzp{\xx}{\vv}$ be the Bongartz completion and $\Bgzm{\xx}{\vv}$ the Bongartz co-completion of a partial cluster $\vv$ with respect to the initial cluster ${\xx}$. 
The flow \eqref{eq:intro-parallel} can be extended to any partial cluster by replacing $\vv$ with $\Bgzpm{\xx}{\vv}$.

Furthermore, the flow is expected to be continuous with respect to $\xevo$. 
For any non-generic point \eqref{eq:intro-generic} with possibly $a_i=0$
that lies in a low dimensional cell, 
we extend the evo flow to $\xevo$ by allowing the coefficients $a_i$ in \eqref{eq:intro-parallel} to take $0$.
Note that the choice of the top cell containing $\xevo$ is not unique.
So the well-definedness of such a limiting approach needs to be proved.
This has been done in \Cref{lem:well-define}.
Here is the final definition.
\begin{definition}
Let $\xevo=\sum_{i=1}^k a_i x_i$ be a point and $\vv$ a partial cluster. The parallel part of the $\xevo$-evo flow $\flowdu{\xevo}$ for $\vv$ is
\begin{equation}\label{eq:intro-parallel3}
\begin{aligned}
\flowdu{\xevo}(\vv)
 =\sum_{i=1}^k a_i\left(
   \normal_\pm(\tgvec_{x_i}^{\Bgzpm{\ol{\xx}}{\vv}})
  -\normal_\mp(\tgvec_{x_i}^{\Bgzpm{\ol{\xx}}{\vv}})\right),
\end{aligned}
\end{equation}
where $\ol{\xx}$ is any cluster containing $\supp\xevo$. 
\end{definition}

In order to obtain a direction field, we make a radial modification. 

\begin{definition}
Let $\xevo$ be a point in $\CCx(\seeds)$. \emph{The $\xevo$-evo flow} $\flowDU{\xevo}$ at $\ccpt{V}$ is
\begin{gather}\label{eq:intro-EVO}
    \flowDU{\xevo}(\ccpt{V})=
    \flowdu{\xevo}(\vv)-\csum(\flowdu{\xevo}(\vv))\ccpt{V},
\end{gather}
where $\csum(\flowdu{\xevo}(\vv))$ is the coefficient sum of $\flowdu{\xevo}(\vv)$, cf. \eqref{eq:csum}.
\end{definition}

The $\xevo$-evo flow induces an \emph{$\xevo$-foliation} on the cluster complex, cf. \Cref{thm:X-foli}.
\begin{theorem}
For any point $\xevo$ in $\CCx(\seeds)$, the $\xevo$-evo flow induces a piecewise linear \emph{$\xevo$-foliation} with only two singularities $\{\xevo,\Xevo\}$ on $\CCx(\seeds)$.
\begin{gather}
    \CCx(\seeds)\setminus\{\xevo,\Xevo\}=\bigsqcup_{L~\text{$\xevo$-leaves}} L.
\end{gather}
\end{theorem}
One can think of the cells as homogeneous mediums and $\xevo$-leaves as rays of light.
When an X-leaf reaches the boundary of a cell, it must pass into an adjacent cell. 
This cell-wall is like the boundary of two mediums with different refractive index, and it is exactly here that the refraction occurs. 

We summarize the correspondence between discrete combinatorics and the continuous geometry in \Cref{table:D-vs-C}.

\begin{table}[H]\label{table:D-vs-C}
\centering
\small
\renewcommand{\arraystretch}{1.2}
\begin{tabular}{@{}p{.45\linewidth}p{.49\linewidth}@{}}
\hline
\textbf{Discrete combinatorics} & \textbf{Continuous geometry}\\
\hline
Cluster exchange graph $\uCEG(\seeds)$
 & Cluster complex $\CCx(\seeds)$\\
Initial cluster $\xx\in\uCEG_0(\seeds)$
 & Source point $\xevo\in\CCx(\seeds)$ (w.r.t. $\flowD{\xevo}$)\\
Green/red mutation
 & Downward/upward cell-crossing\\
$\xx$-orientation $\CEG^{\xx}(\seeds)$
 & Oriented $\xevo$-foliation\\
Source $\xx$ and sink $\DT{\xx}$
 & Singular source $\xevo$ and sink $\Xevo$\\
Bongartz completion/co-completion
 & Upward/downward $\xevo$-support\\
Completion w.r.t a cluster
 & Partial completion w.r.t. $\supp(\xevo)$\\
Wall-crossing (codimension$=1$)
 & Cell-crossing (codimension$\ge1$)\\
 
\hline
\end{tabular}
\caption{From green mutation to $\xevo$-evolution.}
\label{table}
\end{table}

\subsection{Main results}\


In this paper, our main motivation is to prove the following conjecture in the cluster-algebraic setting.
\begin{conjecture}\cite[Conj.~1.7]{QT}\label{Conj0}
    Every $\xevo$-foliation is acyclic.
\end{conjecture}
We say an $\xevo$-foliation is \emph{acyclic} if no $\xevo$-leaf is homotopic to $\sphere^1$.
Under a mild assumption \Cref{ass0}, the acyclicity of $\xevo$-foliation in the cluster-algebraic setting is shown in \Cref{thm:acyclic}.
As an interesting byproduct, we introduce the $\xevo$-supports by the $\xevo$-evo flow.
It turns out to coincide with a partial version of Bongartz (co-)completions defined by intersections, cf. \Cref{thm:X=xx}.
In particular, the rank invariant of any $\xevo$-leaf is well-defined,

\paragraph{\textbf{$\xevo$-supports as partial Bongartz (co-)completions}}\

We introduce the \emph{$\xevo$-support} 
as the support of the $\xevo$-evo flow originated from a partial cluster.

\begin{definition}
    The \emph{upward/downward $\xevo$-support} of a partial cluster $\vv$ is
    \begin{gather}\label{eq:intro-Xsupp}
        \suppdu{\xevo}(\vv)
        =\vv\cup\supp(\flowdu{\xevo}(\vv))
    \end{gather}
    Since the $\xevo$-support depends only on the support of $\xevo$, one can write $\xx$ instead of $\xevo$ in \eqref{eq:intro-Xsupp}.
\end{definition}
For any cluster $\ol{\xx}$, the $\ol{\xx}$-support is precisely the Bongartz (co-)completion with respect to $\ol{\xx}$, i.e.
\begin{gather}\label{eq:intro-supp=Bgz}
    \suppdu{\ol{\xx}}(\vv)=\Bgzpm{\ol{\xx}}{\vv},
\end{gather}
Taking $\xx$ to be a partial cluster, we naturally find a partial variation of Bongartz (co-)completions. 
It is interesting that this partial version can be defined intrinsically as the intersection of larger Bongartz (co-)completions.
\begin{thmdef}\label{thm:intro-supp=Bgz}
    For any partial clusters $\xx,\vv$, the \emph{partial Bongartz (co-)completions} are defined by
    \begin{gather}
        \Bgzpm{\xx}{\vv}:=\bigcap_{\xx\subset\ol{\xx}}\Bgzpm{\ol{\xx}}{\vv},
    \end{gather}
    where $\ol{\xx}$ runs through all clusters containing $\xx$. Then
    \begin{gather}
        \suppdu{\xx}(\vv)=\Bgzpm{\xx}{\vv},
    \end{gather}
    In particular,
    \begin{gather}
        |\suppd{\xx}(\vv)|=|\suppu{\xx}(\vv)|\ge|\xx|.
    \end{gather}
\end{thmdef}

\paragraph{\textbf{Rank invariant}}\

The $\xevo$-foliation decomposes the cluster complex into 1-dimensional piecewise linear manifolds, i.e., the \emph{$\xevo$-leaves}.
The $\xx$-supports illustrate a rank invariant for each $\xevo$-leaf. 
When cell-crossing happens for an oriented $\xevo$-leaf $L$, the downward flow passes from $\cc(\suppu{\xevo}(\vv))$ to $\cc(\suppd{\xevo}(\vv))$. 
\Cref{thm:intro-supp=Bgz} shows that the two cells have the same rank, even when the crossing wall has higher codimension. 
Thus we define the following leaf rank.
\begin{cordef}
    Let $L$ be an $\xevo$-leaf and $\cc(\vv)^\circ$ intersects $L$ nontrivially. 
    The \emph{rank} of $L$ is
    \begin{gather}
        \rank L:=|\suppu{\xevo}(\vv)|=|\suppd{\xevo}(\vv)|,
    \end{gather}
    which does not depend on the choice of $\cc(\vv)$ intersecting $L$ nontrivially.
\end{cordef}
For generic $\xevo$, every $\xevo$-leaf has rank $n$. 
For non-generic $\xevo$, there are lower-rank $\xevo$-leaves living in lower dimensional cells. They are not visible in the exchange graph alone.

\paragraph{\textbf{Acyclicity}}\



Using Cao's dominant sets, cf. \cite{Cao}, we prove the main theorem of this paper.
\begin{theorem}
    Assuming \Cref{ass0},
    every $\xevo$-foliation is acyclic.
\end{theorem}

\subsection{Organization of the paper}\label{sec:further}\


This paper is organized as follows.
\begin{itemize}
    \item \Cref{sec:preCA} recalls the cluster-algebraic preliminaries.
    \item \Cref{sec:suppBgz} constructs the parallel fields and identifies $\xevo$-supports with partial Bongartz (co-)completions.
    \item \Cref{sec:evoCC} introduces the radial correction, constructs the foliation and establishes the cell-crossing and rank properties.
    \item \Cref{sec:acyclic} develops the order-theoretic approach to acyclicity and its topological consequences.
    \item Finally, \Cref{app:cats} explains the categorical construction and its relationship with $\tau$-tilting theory.
\end{itemize}

\subsection{Notation list}


\begingroup
\renewcommand{\arraystretch}{1.15}
\begin{longtable}{@{}p{.32\linewidth}p{.64\linewidth}@{}}
\hline
Notation & Meaning \\
\hline\endhead
$\seeds$& A cluster pattern. \\
$\uCEG(\seeds)$ & The cluster exchange graph. \\
$\CEG^{\xx}(\seeds)$& The $\xx$-orientation of $\uCEG(\seeds)$. \\
$\CCx(\seeds)$ & The cluster complex. \\
$\CCx_*(\seeds)$ & The set of partial clusters in $\seeds$. \\
$\cc(\vv)$ & The cell of $\vv$. \\
$\cc(\vv)^\circ$ & The interior of $\cc(\vv)$. \\
$\dcc(\vv)$ & The tangent directions of $\cc(\vv)$. \\
$\supp(\ccpt{V})$; $\supp(\cvec)$ & The support of a point $\ccpt{V}$ or a vector $\cvec$. \\
$\normal_\pm(\cvec)$& The normalized positive/negative part of a vector $\cvec$. \\
$\Cmat^{\xx}_{\vv}$; $\Gmat^{\xx}_\vv$ & The $C$- and $G$-matrices with initial cluster $\xx$. \\
$\cvec^{\xx}_{v;\vv}$; $\gvec^{\xx}_v$ & The $c$- and $g$-vectors with initial cluster $\xx$. \\
$\tgvec^{\vv}_x$ & The row of the $C$-matrix $\Cmat^{\xx}_\vv$. \\
$\DT{(-)}$ & The cluster DT-transformation. \\
$\Bgzpm{\xx}{\vv}$& The (partial) Bongartz completion/co-completion of $\vv$ with respect to $\xx$. \\
$\flowdu{\xevo}$ & The upward/downward parallel parts of the $\xevo$-evo flow. \\
$\flowDU{\xevo}$ & The upward/downward $\xevo$-evo flow. \\
$\suppdu{\xevo}(\vv)$ & The upward/downward $\xevo$-support.\\
$\dom{\xx}(\vv)$ & The dominant set of $\vv$. \\
$\le_{\xx}$ & The $\xx$-order.\\
\hline
\end{longtable}
\endgroup

\subsection*{Acknowledgments}
We would like to thank Peigen Cao and Changjian Fu for inspiring conversations.
\paragraph{\textbf{AI disclosure}}
The use of AI is purely cosmetic:
the first draft was finished in late 2025. 
Due to various reasons, we only pick it up until Sep. 2026 where we use AI to improve the writing and presentation.

\section{Preliminaries on cluster algebras}\label{sec:preCA}

\begin{convention}\label{conv:eE}
We adopt the following typographical conventions:
\begin{itemize}
    \item 
    boldface for an $n$-tuple of variables, e.g. $\xx=(x_1,\ldots,x_n)$;
    \item slanted boldface for a column vector, e.g. $\bfit{e}_i$ denotes the $i$th standard column vector;
    \item underlined boldface for a row vector, e.g. $\ubf{e}_i$ denotes the $i$th standard row vector; and
    \item 
    capital letters for matrices, e.g. $E_{ij}=\bfit{e}_i\cdot\ubf{e}_j$ denotes the elementary matrix.
\end{itemize}
Moreover, we write
\begin{itemize}
    \item $(-)^T$ the matrix transpose,
    \item $\bfit{I}$ the identity matrix,
    \item $J_i=\bfit{I} - 2E_{ii}$ the $i$-th reflection matrix.
    \item $\xx^\bfit{h}=\prod_{i=1}^n x_i^{h_i}$ for an $n$-tuple $\xx=(x_1,\ldots,x_n)$ and a column vector $\bfit{h}=(h_1,\ldots,h_n)^T\in\ZZ^{n}$
\end{itemize}
\end{convention}

\subsection{Cluster algebras}\label{sec:CA}\


A \emph{semifield} is an abelian group $(\PP,\cdot)$ endowed with an auxiliary addition $\oplus$ that is distributive over the multiplication $\cdot$.
The only torsion element of $\PP$ is $1$. Consequently, the group ring $\ZZ\PP$ is a domain, and its field of fractions $\QQ\PP$ is well defined.



Let $n$ be a positive integer, and let $\hua{F}$ be a field of rational functions in $n$ independent variables over $\QQ\PP$.
A \emph{(labeled) seed over $\PP$} is a triple $\seed=(\xx,\yy,B)$, where
\begin{itemize}
    \item the \emph{cluster} $\xx=(x_1,\ldots,x_n)$ is a transcendence basis of $\hua{F}$, and a \emph{partial cluster} is a subset of a cluster;
    \item the \emph{coefficient tuple} $\yy=(y_1,\ldots,y_n)$ is an $n$-tuple in $\PP$ and
    \item the \emph{exchange matrix} $B=(b_{ij})$ is an $n\times n$ skew-symmetrizable integer matrix; that is, there exists $D=\on{diag}(d_1,\ldots,d_n)$ with $d_i\in\QQ_{>0}$ such that $DB$ is skew-symmetric.
\end{itemize}
The \emph{$\hat{y}$-variables} of a seed $\seed$ form the $n$-tuple $\hat{\yy}=(\hat{y}_1,\ldots,\hat{y}_n)$ in $\hua{F}$ defined by
\begin{gather}\label{eq:haty}
    \hat{y}_i=y_i\xx^{B\bfit{e}_i}.
\end{gather}

The \emph{mutation of $\seed$ in direction $k$}, denoted by $\mu_k(\seed)=(\xx',\yy',B')$, is defined as follows:
\begin{align}
&x'_i=
\begin{cases}\label{eq:mutx}
x_k^{-1} \Big( \prod_{j=1}^n x_j^{\pos{-b_{jk}}} \Big) \displaystyle{\frac{1 + \hat{y}_k}{1 \oplus y_k}} & i=k,\\
x_i & i\ne k,
\end{cases}
\\
&y'_i=
\begin{cases}
y_k^{-1} & i=k,\\
y_i y_k^{\pos{b_{ki}}}(1 \oplus y_k)^{-b_{ki}} & i\ne k,
\end{cases}
\\
&b'_{ij}=
\begin{cases}
-b_{ij} & i=k~\text{or}~j=k,\\
b_{ij}+b_{ik}\pos{b_{kj}}+\pos{-b_{ik}}b_{kj} & i,j\ne k,
\end{cases}
\end{align}
where $\pos{a}=\max(0,a)$ and $\nega{a}=a-\pos{a}$.
The following standard properties are immediate from the definition.
\begin{itemize}
    \item $\xx'$ is again a transcendence basis of $\hua{F}$; hence the mutation of a seed is again a seed.
    \item The mutation $\mu_k$ is an involution, i.e., $\mu_k^2(\seed)=\seed$ for any seed $\seed$.
    \item If $D$ skew-symmetrizes $B$, then it also skew-symmetrizes $B'$.
\end{itemize}

Let $\regT_n$ be the \emph{$n$-regular tree}, that is, the tree in which every vertex has valency $n$. Label the $n$ edges incident with each vertex by $1,\ldots,n$, with no repetitions. Two vertices $t$ and $t'$ are said to be \emph{$k$-adjacent} if they are joined by an edge labeled $k$.

A \emph{cluster pattern} is a family of rank-$n$ seeds with coefficients
\begin{gather}\label{eq:cpattern}
    \seeds=\{ \seed_t=(\xx_t,\yy_t,B_t) \}_{ t\in\regT_n }
\end{gather}
indexed by $\regT_n$ and satisfying $\seed_t=\mu_k(\seed_{t'})$ whenever $t$ and $t'$ are $k$-adjacent.
A skew-symmetrizer $D$ of $\seeds$ is the common skew-symmetrizer of all seeds in $\seeds$.
By \cite[Prop.~3]{CL2}, every seed in a cluster pattern $\seeds$ is uniquely determined by its cluster. We may therefore write $\xx$ for the seed $(\xx,\yy,B_\xx)$ in a fixed cluster pattern.

\begin{definition}\label{def:CEG}
The \emph{(unoriented) cluster exchange graph} $\uCEG(\seeds)$ of a cluster pattern $\seeds$ is an $n$-regular graph whose vertices are clusters (up to permutation) and edges are mutations.
\end{definition}

The \emph{reduced cluster pattern} with respect to a partial cluster $\vv$ is
\begin{gather}\label{eq:red-pattern}
    \seeds_\vv=\{ \seed_t=(\xx,\yy,B)    \mid    \vv\subset\xx \},
\end{gather}
where the seeds mutate only at unfrozen variables in $\xx\backslash\vv$.
Its exchange graph $\uCEG(\seeds_\vv)=\uCEG(\seeds,\vv)$ is the full subgraph of $\uCEG(\seeds)$ consisting of the clusters containing $\vv$, which is connected by \cite[Thm.~10]{CL2}.


\subsection{$C$-Matrices, $G$-matrices and $g$-fans}\label{sec:gc}\

Let $\seeds$ be a cluster pattern, let $\xx=(\xx,\yy,B_\xx)$ be its initial seed, and let $D$ be a skew-symmetrizer of $\seeds$.
An important fact is the Laurent-positivity phenomenon: after using the separation formulas, each cluster variable is a Laurent polynomial in the initial cluster variables with nonnegative coefficients (in the appropriate coefficient ring). 

Moreover, the cluster pattern $\seeds=\{(\vv,\zz,B_{\vv})\}$ can be reconstructed from the following inductively defined invariants:
\begin{itemize}
    \item the \emph{$C$-matrices} $\Cmat^\xx_\vv$, which encode the tropical part of $\zz$;
    \item the \emph{$G$-matrices} $\Gmat^\xx_\vv$, which encode the tropical part of $\vv$; and
    \item the \emph{$F$-polynomials} $\Fpoly^\xx_\vv$, which encode the nontropical parts of both $\vv$ and $\zz$.
\end{itemize}
We first recall the definition of $C$/$G$-matrices.

\begin{definition}\label{def:CG}
Let $\Cmat^{\xx}_{\vv}$/$\Gmat^{\xx}_{\vv}$ be the \emph{$C$-matrices/$G$-matrices} of $\vv$ with respect to $\xx$.
Their rows and columns are indexed by variables of $\xx=\{x_i\}_{i=1}^n$ and $\vv=\{v_j\}_{j=1}^n$, respectively.
These matrices are uniquely determined by
\begin{itemize}
  \item the initial condition that $\Cmat^\xx_\xx=\bfit{I}=\Gmat^\xx_\xx$ and
  \item the inductive mutation formula for $\vv'=\mu_{v_k}(\vv)$
    \begin{gather}
    \label{eq:mutC} \Cmat^{\xx}_{\vv'}=\Cmat^{\xx}_{\vv}J_{v_k}+\Cmat^{\xx}_{\vv}E_{v_k,v_k}\pos{B_\vv}+\pos{-\Cmat^{\xx}_{\vv}}E_{v_k,v_k} B_\vv.\\
    \label{eq:mutG}\Gmat^{\xx}_{\vv'}=\Gmat^{\xx}_{\vv}J_{v_k}+\Gmat^{\xx}_{\vv}\pos{-B_\vv}E_{v_k,v_k}- B_{\xx}\pos{-\Cmat^{\xx}_{\vv}}E_{v_k,v_k}.
    \end{gather}
The indexing of the variables of $\vv'=\mu_k(\vv)$ is inherited from that of $\vv$.
\end{itemize}
The column vectors of a $C$-matrix/$G$-matrix are called the \emph{$c$-vectors/$g$-vectors}, denoted by $\cvec^\xx_{v_i,\vv}$/$\gvec^\xx_{v_i,\vv}$, respectively.
It is known that the $g$-vectors do not depend on the choice of $\vv$ so we can write $\gvec^\xx_{v_i}=\gvec^\xx_{v_i,\vv}$,
cf. \cite[Part.~I, Def.~4.8]{Nak}.

For convenience, we will define the following transposed $C/G$-matrices.
\begin{gather}\label{eq:dual3}
    \widetilde{\Cmat}^{\vv}_{\xx}\colon=(\Gmat^{\xx}_\vv)^T
    \quad\text{and}\quad
    \widetilde{\Gmat}^{\vv}_{\xx}\colon=(\Cmat^{\xx}_\vv)^T. 
\end{gather}
They are indeed the $C/G$-matrices for the cluster pattern generated by the initial seed $(\vv,B^T_\vv)$, which shares the same cluster exchange graph with $\seeds$, cf. \cite{NZ}.
Thus the notation is well-defined.
Their column vector, 
denoted by $\widetilde{\cvec}^{\vv}_{x,\xx}$/$\tgvec^{\vv}_x$, are precisely the row of $\Gmat^\xx_\vv$/$\Cmat^\xx_\vv$.
\end{definition}

The columns of the $C$- and $G$-matrices form $\ZZ$-bases; see \cite[Prop.~4.2]{NZ}.
\begin{gather}\label{eq:det1}
    \det(\Cmat^\xx_\vv)=\pm1=\det(\Gmat^\xx_\vv).
\end{gather}

We recall several tropical dualities between $C$- and $G$-matrices; see \cite{FZ,NZ}.

\begin{lemma}
We have the following dualities.
\begin{align}
    \label{eq:dual1}
    \Gmat^{\xx}_\vv B_{\vv}&= B_{\xx} \Cmat^{\xx}_\vv,\\
    \label{eq:dual2}
    (\Gmat^{\xx}_\vv)^T D \Cmat^{\xx}_\vv=&~D=(\Cmat^{\xx}_\vv)^T D \Gmat^{\xx}_\vv.
\end{align}
\end{lemma}


The geometry of the $G$-matrices is captured by $g$-cones and $g$-fans.
After multiplying $D$ by a common positive integer if necessary, we assume that $D$ has positive integral diagonal entries. 
Let $M\cong\ZZ^n$ and $N\cong\ZZ^n$, equipped with the bilinear pairing
\begin{gather}\label{eq:pairing}
    (\uu,\ww)_D=\uu^T D \ww.
\end{gather}
The pairing extends $\RR$-bilinearly to $M_\RR\times N_\RR$, where $M_\RR=M\otimes_\ZZ\RR$ and $N_\RR=N\otimes_\ZZ\RR$.
For vectors $\bfit{a}_1,\ldots,\bfit{a}_k\in M_\RR$, we denote the convex cone that they generate by
\[
    \con(\bfit{a}_1,\ldots,\bfit{a}_k)=\RR_{\ge0} \bfit{a}_1 + \ldots + \RR_{\ge0} \bfit{a}_k.
\]
\begin{definition}
The \emph{$g$-cone} associated with a $G$-matrix $\Gmat^\xx_\vv$ is the convex cone in $M_\RR$ generated by its $g$-vectors:
\begin{gather}\label{eq:gcone}
    \con(\Gmat^\xx_\vv)=\con(\gvec^\xx_{v_1},\ldots,\gvec^\xx_{v_n}).
\end{gather}
The \emph{$g$-fan} of a cluster pattern $\seeds$ with respect to an initial seed $\xx$ is the collection of convex cones generated by compatible $g$-vectors:
\begin{gather}\label{eq:gfan}
    \fan^\xx(\seeds)=\{\con(\gvec^\xx_{v_1},\ldots,\gvec^\xx_{v_k}):\{v_1,\ldots,v_k\}\text{ is a partial cluster}\}.
\end{gather}
Its maximal cones are precisely the $g$-cones. We denote the ground set of $\fan^\xx(\seeds)$ by $|\fan^\xx(\seeds)|$.
\end{definition}

\begin{remark}\label{rem:fan-CSD}
For our conventions concerning cones and fans, see \cite[Part~II, Sec.~2]{Nak}. We will use the following facts.
\begin{enumerate}
    \item By the tropical duality \eqref{eq:dual2}, the codimension-one cones of the $g$-fan have $c$-vectors as inward-pointing normal vectors with respect to the pairing \eqref{eq:pairing}.
    \item The $g$-fans of $\seeds$ with initial seed $\xx$ are associated with the cluster scattering diagram of $\seeds$ with initial seed $\xx$ as the reachable chambers, cf. \cite[Lem.~2.10]{GHKK18}.
\end{enumerate}
\end{remark}

\subsection{Green mutation}\

We retain the notation of \Cref{sec:CA}.
Let $\ZZp$ and $\ZZm$ denote the sets of nonnegative and nonpositive integers, respectively.
Sign-coherence, conjectured in \cite[Conj.~6.13]{FZ07} and proved in \cite{GHKK18}, is one of the fundamental properties of cluster algebras.

\begin{lemma}\cite[Cor.~5.5]{GHKK18}\label{lem:sc}
The $C$-matrices are \emph{column sign-coherent}, that is,
\begin{gather}\label{eq:c-sc}
\cvec^{\xx}_{v;\vv}\in\ZZp^n~\text{or}~\ZZm^n.
\end{gather}
Dually, the $G$-matrices are \emph{row sign-coherent}, that is,
\begin{gather}\label{eq:g-sc}
    \widetilde{\cvec}^\vv_{x,\xx}=\ubf{e}_{x} \Gmat^\xx_\vv\in\ZZp^n~\text{or}~\ZZm^n.
\end{gather}
\end{lemma}

Sign-coherence induces an orientation on $\uCEG(\seeds)$, known as the \emph{green-mutation orientation}; see \cite[\S5]{Ke2}. 

\begin{definition}\label{def:green}
A cluster variable $v\in\vv$ is \emph{green} (respectively, \emph{red}) with respect to $\xx$ if $\cvec^{\xx}_{v;\vv}\in\ZZp^n$ (respectively, $\cvec^{\xx}_{v;\vv}\in\ZZm^n$).
The \emph{$\xx$-orientation} $\CEG^\xx(\seeds)$ of $\uCEG(\seeds)$ is obtained by orienting each edge in the direction of mutation at a green cluster variable.
Note that
\begin{itemize}
    \item {$\vv$ is a source of $\CEG^\xx(\seeds)$ if and only if every $v\in\vv$ is green;}
    \item {$\vv$ is a sink of $\CEG^\xx(\seeds)$ if and only if every $v\in\vv$ is red.}
\end{itemize}
\end{definition}

The following notion is useful in discussing the existence and uniqueness of the source and sink of $\CEG^\xx(\seeds)$.

\begin{definition}\cite[Def.~3.1.1]{Mul}.
A \emph{green-to-red sequence} (or \emph{reddening sequence}) with respect to $\xx$ is a mutation sequence starting at $\xx$ and ending at a cluster $\xx'$ all of whose cluster variables are red with respect to $\xx$. 
The individual mutations in the sequence need not all be green.
\end{definition}

The exchange graph is the adjacency graph of the maximal 
cones of $\fan^\xx(\seeds)$, where two maximal cones are adjacent precisely when they share a facet. This gives the following properties.

\begin{proposition}\label{prop:ss}
For each cluster $\xx$,
\begin{itemize}
    \item $\CEG^\xx(\seeds)$ admits a unique source $\xx$.
    \item $\CEG^\xx(\seeds)$ admits a unique sink $\Vec{\mu}(\xx)$ if and only if there exists a green-to-red sequence $\Vec{\mu}$ with respect to $\xx$.
\end{itemize}
In particular, we have
\begin{gather}\label{eq:C-C-dual0}
    \Cmat^\xx_{\Vec{\mu}(\xx)}=-\bfit{P}
\end{gather}
for some permutation matrix $\bfit{P}$.
Thus $\xx$ and $\Vec{\mu}(\xx)$ may be canonically identified via
\begin{gather}\label{eq:id-CC}
  x(\in\xx)
  \xlongleftrightarrow{\cvec^\xx_{x,\xx}=-\cvec^\xx_{x',\Vec{\mu}(\xx)}}
  x'(\in\Vec{\mu}(\xx)).
\end{gather}
\end{proposition}


The conjugation of a green-to-red sequence by mutations is another green-to-red sequence, cf. \cite[Thm~3.2.1]{Mul}.
If $\seeds$ admits a green-to-red sequence with respect to a cluster $\xx$, then it admits a green-to-red sequence with respect to any other cluster.
Under \Cref{ass0}, one can define the \emph{cluster DT-transformation}; see \cite{Nag}.

\begin{definition}\label{def:DT}
Denote by $\Vec{\mu}^\xx$ the green-to-red sequence with respect to $\xx$ in $\seeds$.
Then the cluster DT-transformation of a cluster is defined by
\begin{gather}
    \DT{(-)}:\seeds\to\seeds
\end{gather}
sending each cluster $\xx$ to $\Vec{\mu}^\xx(\xx)$.
Since the conjugation of $\Vec{\mu}^\xx$ by a mutation is another green-to-red sequence.
The operation $\DT{(-)}$ commutes with mutations and can be specified on cluster variables as follows.
\begin{gather}\label{eq:cDT}
    \DT{(-)}:\{\text{cluster variables in }~\seeds\}\to\{\text{cluster variables in }~\seeds\}
\end{gather}
\end{definition}


\subsection{$F$-polynomials}\label{sec:F-poly}\

We next recall the definition of $F$-polynomials, the second key ingredient in the separation formula for cluster variables.

\begin{definition}\label{def:Fpoly}
Let $\Fpoly^\xx_{\vv}(\yy)=(\Fpoly^\xx_{v_1}(\yy),\ldots,\Fpoly^\xx_{v_n}(\yy))\in\ZZ[\yy]^{n}$ be the $F$-polynomials for a cluster $\vv$ with respect to $\xx$.
They are uniquely determined by
\begin{itemize}
  \item the initial condition $\Fpoly^\xx_x(\yy)=1$ for any $x\in\xx$ and
  \item the inductive mutation formula for $\vv'=\mu_{v_k}(\vv)=\{v_i\}_{i\ne k}\cup\{v'_k\}$
    \begin{gather}\label{eq:mutF}
        \Fpoly^\xx_{v'_k}\Fpoly^\xx_{v_k}=
        \yy^{[\cvec^{\xx}_{v_k,\vv}]_+} \Fpoly^{\xx}_{\vv}(\yy)^{[B_{\vv}\bfit{e}_{v_k}]_+} +
        \yy^{[-\cvec^{\xx}_{v_k,\vv}]_+} \Fpoly^{\xx}_{\vv}(\yy)^{[-B_{\vv}\bfit{e}_{v_k}]_+}.
    \end{gather}
\end{itemize}
\end{definition}

We will use the following standard properties of $F$-polynomials; see \cite[Props.~3.6 and~5.2]{FZ07} and \cite[Cor.~0.4]{GHKK18}:
\begin{itemize}
    \item $\Fpoly^\xx_v(\yy)\in\ZZ[\yy]$,
    \item $\Fpoly^\xx_v(0)=1$,
    \item all coefficients of $\Fpoly^\xx_v(\yy)$ are nonnegative.
\end{itemize}

We now recall the \emph{separation formula}.

\begin{theorem}\label{thm:CGF}\cite[\S3]{FZ07}
Let $\seeds$ be a cluster pattern, and let $\xx=(\xx,\yy,B_\xx)$ be an initial seed with $\hat{y}$-variables $\hat{\yy}$. Then, for any seed $\vv=(\vv,\zz,B_\vv)$ in $\seeds$,
\begin{gather}
    v_i=\xx^{\gvec^\xx_{v_i}}\cdot\frac{\Fpoly^\xx_{v_i}(\hat{\yy})}{\Fpoly^\xx_{v_i}|_{\PP}(\yy)}
    \quad\text{and}\quad
    z_i=\yy^{\cvec^\xx_{v_i,\vv}}\cdot\Fpoly^\xx_\vv|_{\PP}(\yy)^{B_\vv\bfit{e}_i}
\end{gather}
\end{theorem}

\subsection{Cluster complex}\



We study the cluster complex of $\seeds$ as a combinatorial counterpart of its $g$-fans.

\begin{definition}\label{def:Cpx}
The \emph{cluster complex} $\CCx(\seeds)$ of a cluster pattern $\seeds$ is the clique complex associated with the compatibility relation on the set of cluster variables. Its simplices are the partial clusters, and its maximal simplices are the clusters.

We denote by $\CCx_k(\seeds)$ the set of partial clusters of rank $k$, and by $\CCx_*(\seeds)=\bigsqcup\CCx_k(\seeds)$ the set of all partial clusters.
\end{definition}

For any choice of initial cluster $\xx$, the cluster complex $\CCx(\seeds)$ admits a natural embedding into the $g$-fan $\fan^\xx(\seeds)$:
\begin{equation}\label{eq:gf-embedding}
\begin{array}{rccc}
  \gvec^\pcts{\xx}\colon
    & \CCx(\seeds) & \longrightarrow & |\fan^{{\xx}}(\seeds)|  \\
    & \sum_{i=1}^k c_i v_i & \longmapsto & \sum_{i=1}^k c_i \gvec^\xx_{v_i}.
\end{array}
\end{equation}
In particular, this gives the following one-to-one correspondence between the partial clusters in $\CCx(\seeds)$ and the cones in $\fan^\xx(\seeds)$:
\begin{equation}
\begin{array}{rccc}
  \sigma^\pcts{\xx}\colon
    & \CCx_*(\seeds) & \longrightarrow & \fan^{{\xx}}(\seeds)  \\
    & \vv & \longmapsto & \sigma( \gvec^\xx_{v} )_{v\in\vv}.
\end{array}
\end{equation}

\begin{figure}
    \centering
    \reviewfigure[width=\linewidth]{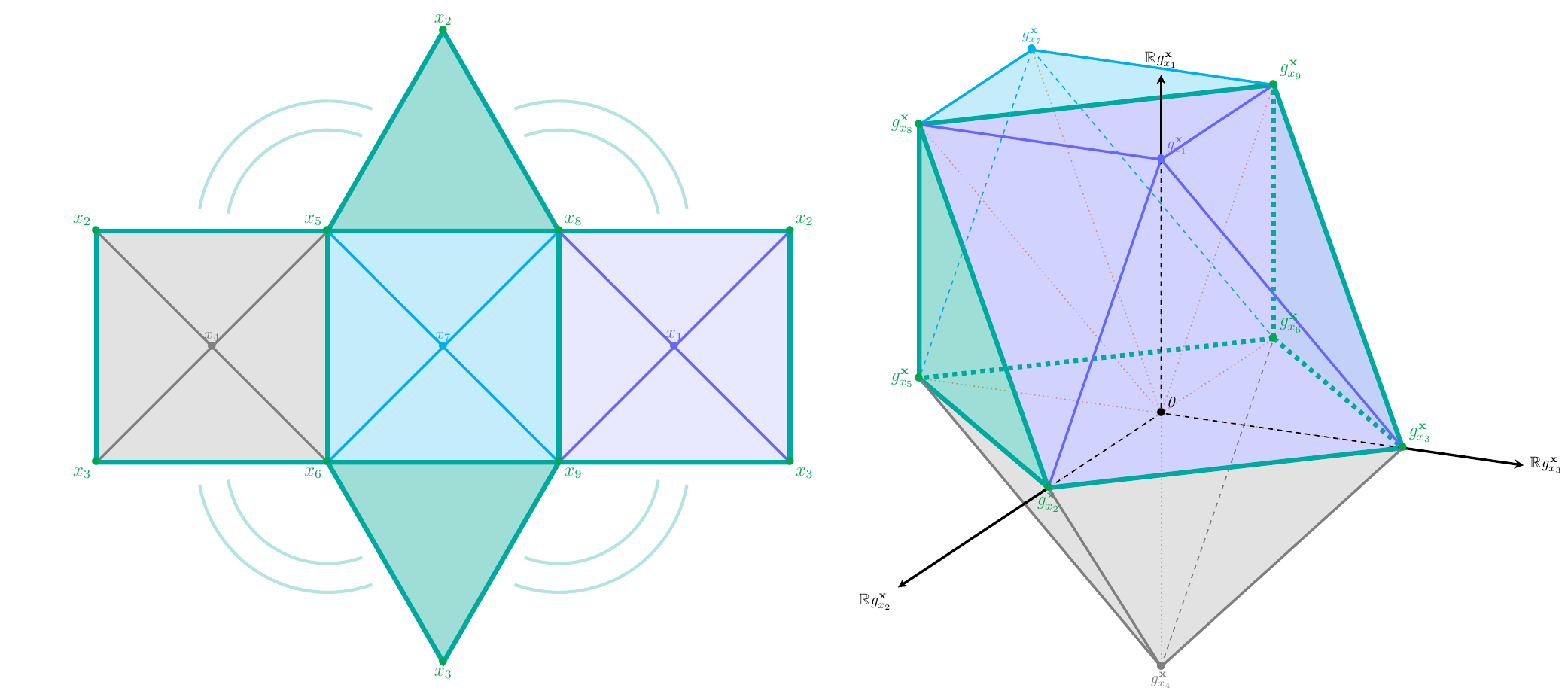}
    \caption{The cluster complex and $g$-fan of $\seeds_{\DynA{3}}$}
    \label{fig:A3-g-fan}
\end{figure}

For example, \Cref{fig:A3-g-fan} illustrates the cluster complex $\CCx(\DynA{3})=\CCx(\seeds_{\DynA{3}})$ and its natural embedding into the $g$-fan $\fan^\xx(\seeds_{\DynA{3}})$. 

For each partial cluster $\vv\in\CCx_*(\seeds)$, we use the following notations.
\begin{itemize}
    \item The formal sum with respect to $\vv=(v_1,\ldots,v_k)$ is 
    \begin{gather}
        \RR\vv=\{   \cvec=\sum_{i=1}^k c_i v_i    \mid    c_i\in\RR\}.
    \end{gather}
    It is similar for $\RR_{\ge0}\vv$ and $\RR_{>0}\vv$.
    The formal sum can be viewed as a column vector indexed by $\vv$, denoted by $\cvec=\sum_{i=1}^k c_i v_i$.
    For example, $\gvec^\vv_{x_i}\in\RR\vv$ when $\vv$ is a cluster.
    
    
    \item The coefficient sum of a formal sum is 
    \begin{equation}\label{eq:csum}
    \begin{array}{rccl}
      \csum\colon
    & \RR\vv & \longrightarrow & \RR  \\
    & \cvec & \longmapsto & \csum(\cvec)=\sum_{i=1}^k c_i.
    \end{array}
    \end{equation}

    \item The positive and negative parts of $\cvec=\sum_{i=1}^k c_i v_i$ is defined by
    \begin{gather}\label{eq:pmvec}
        [\cvec]_\pm=\sum_{i=1}^k [c_i]_\pm v_i.
    \end{gather}
\end{itemize}
In the cluster complex $\CCx(\seeds)$,
\begin{itemize}
    \item The \emph{cell} $\cc(\vv)$ associated with $\vv$ is the simplex
    \begin{gather}\label{eq:cc}
    \cc(\vv)=\{\cvec\in\RR_{\ge0}\vv \mid \csum(\cvec)=1 \}.
    \end{gather}
    Its relative \emph{interior} $\cc(\vv)^\circ=\cc(\vv)\cap\RR_{>0}\vv$,
    which agrees with the usual notion of relative interior except when $k=1$, we set $\{v_1\}^\circ=\{v_1\}$. 
    A \emph{generic point} is an interior point of a top-dimensional cell.

     \item
    The space of \emph{(tangent) directions} of the cell $\cc(\vv)$ is
    \begin{gather}\label{eq:dcc}
        \dcc(\vv)=\{\cvec\in\RR\vv \mid \epsilon(\cvec)=0\}.
    \end{gather}
    The tangent space at any point of $\cc(\vv)^\circ$ is naturally identified with $\dcc(\vv)$.
    
\end{itemize}
For a cluster complex $\CCx(\seeds)$, we denote by
\begin{gather}\label{eq:direct}
    \sCpx(\seeds)=\bigcup_{\vv\in\CCx_*(\seeds)}\RR\vv
    \quad\text{and}\quad
    \dCpx(\seeds)=\bigcup_{\vv\in\CCx_*(\seeds)}\dcc(\vv).
\end{gather}

Note that $\csum$ and $[-]_\pm$ can be naturally defined on $\sCpx(\seeds)$. 
We moreover need the following operators for the vectors in $\sCpx(\seeds)$:
\begin{itemize}
    \item The \emph{support} of a vector is 
    \begin{gather}\label{eq:supp}
    \begin{array}{rccc}
    \supp\colon
    & \sCpx(\seeds) & \longrightarrow & \CCx_*(\seeds)  \\
    & \cvec=\sum_{i=1}^k c_i v_i & \longmapsto & [v_i]_{c_i\ne0},
    \end{array}
    \end{gather}
    In particular, $\supp(\ccpt{V})=\vv$ if and only if $\ccpt{V}\in\cc(\vv)^\circ$.

    \item The \emph{extended DT-transformation} (cf. \eqref{eq:cDT}) is 
    \begin{gather}\label{eq:extendDT}
    \begin{array}{rccc}
    \DT{(-)}\colon
    & \sCpx(\seeds) & \longrightarrow & \sCpx(\seeds)  \\
    & \cvec=\sum_{i=1}^k c_i v_i & \longmapsto & \DT{\cvec}=\sum_{i=1}^k c_i \DT{v_i},
    \end{array}
    \end{gather}

    \item The \emph{normalized positive/negative part} is an operator
    \begin{gather}\label{eq:op-normal}
    \begin{array}{rccc}
    \normal_\pm\colon
    & \sCpx(\seeds) & \longrightarrow & \CCx(\seeds)\cup\{\bfit{0}\}
    \end{array}
    \end{gather}
    sending a vector $\cvec\in\sCpx(\seeds)$ to
    \begin{equation}\label{eq:normal}
        \normal_\pm(\cvec)=
        \begin{cases}
            [\cvec]_\pm/\csum([\cvec]_\pm),
            &\csum([\cvec]_\pm)\ne0\\
            \bfit{0},
            &\csum([\cvec]_\pm)=0.
        \end{cases}
    \end{equation}
\end{itemize}






\section{\texorpdfstring{$\xevo$-supports as partial Bongartz (co-)completions}{X-supports and partial Bongartz (co-)completions}}\label{sec:suppBgz}

In this section, we will show the local  $\xevo$-foliations can be glued via extending the definition of the $\xevo$-evo flow. 
To do this, we generalize the Bongartz (co-)completions introduced by \cite{CGY} to a partial version.

\subsection{Partial Bongartz (co-)completions}\

We first recall the Bongartz (co-)completions for cluster algebras.

\begin{definition}\cite[Def.~4.9]{CGY}\label{def:bong}
Let $\xx$ be a cluster and $\vv$ a partial cluster in $\seeds$.
\begin{itemize}
    \item A cluster $\Bgzp{\xx}{\vv}$ is the \emph{Bongartz completion} for $\vv$ with respect to $\xx$ if
    \begin{equation}\label{eq:bong}
    \vv\subset\Bgzp{\xx}{\vv}\quad\text{and}\quad
    \cvec^\xx_{v,\Bgzp{\xx}{\vv}}\in\ZZp^n,
    ~\forall v\in\Bgzp{\xx}{\vv}\backslash\vv.
    \end{equation}
    \item A cluster $\Bgzm{\xx}{\vv}$ is the \emph{Bongartz co-completion} for $\vv$ with respect to $\xx$ if
    \begin{equation}\label{eq:cbong}
    \vv\subset\Bgzm{\xx}{\vv}\quad\text{and}\quad
    \cvec^\xx_{v,\Bgzm{\xx}{\vv}}\in\ZZm^n,
    ~\forall v\in\Bgzm{\xx}{\vv}\backslash\vv.
\end{equation}
\end{itemize}
We will also use the \emph{transposed Bongartz (co-)completion} $\tBgzpm{\xx}{\vv}$.
It replaces the $c$-vectors in \eqref{eq:bong}/\eqref{eq:cbong} by the transposed $c$-vectors, cf. \eqref{eq:dual3}.
\end{definition}

For example, we have $\Bgzpm{{\xx}}{\vv}=\vv$ when $\vv$ is a cluster.
In general, the existence and uniqueness of $\Bgzpm{{\xx}}{\vv}$ are associated with the reduced cluster pattern, cf. \eqref{eq:red-pattern}.

\begin{lemma}\label{lem:Bong}\cite[Thm.~5.4, Cor.~5.5]{CGY}
Let $\CEG^\xx(\seeds,\vv)\subset\CEG^\xx(\seeds)$.
Then there exists a graph isomorphism
\begin{gather}
    \CEG^\xx(\seeds,\vv)\cong\CEG^{\Bgzp{\xx}{\vv}\backslash\vv}(\seeds_\vv)
\end{gather}
sending $\uu\in\uCEG(\seeds,\vv)_0$ to $\uu\backslash\vv\in\uCEG(\seeds_\vv)_0$.
Moreover,
    \begin{itemize}
    \item $\Bgzp{\xx}{\vv}$ is the unique source of the induced orientation $\CEG^\xx(\seeds,\vv)$
    (which always exists).
    \item $\Bgzm{\xx}{\vv}$ is the unique sink of the induced orientation $\CEG^\xx(\seeds,\vv)$
    provided its existence for $\vv$ if and only if $\CEG^{\Bgzp{\xx}{\vv}\backslash\vv}(\seeds_\vv)$ admits a sink.
\end{itemize}
\end{lemma}

The existence condition of $\Bgzm{\xx}{\vv}$ for any $\vv$ can be ensured by \Cref{ass0}, cf. \cite[Thm.~3.1.3]{Mul}.

We need the following technical lemma later, which says that when fixing $\vv$ and mutating $\xx$ once,
the corresponding Bongartz (co-)completions differ by at most one mutation.

\begin{lemma}\cite[Lem.~3.10]{FGL}\label{lem:Bgz-mut}
Let $\xx$ be a cluster and $\vv$ a partial cluster.
For any $x\in\xx$, we have the following two cases.
\begin{itemize}
    \item If $\forall v\in\Bgzpm{\xx}{\vv}\backslash\vv,~\cvec^\xx_{v,\Bgzpm{\xx}{\vv}}\ne\pm\bfit{e}_x$, then $\Bgzpm{\xx}{\vv}=\Bgzpm{\mu_x(\xx)}{\vv}$.
    \item If $\exists v\in\Bgzpm{\xx}{\vv}\backslash\vv~\text{s.t.}~\cvec^\xx_{v,\Bgzpm{\xx}{\vv}}=\pm\bfit{e}_x$, then $\mu_{v}(\Bgzpm{\xx}{\vv})=\Bgzpm{\mu_x(\xx)}{\vv}$.
\end{itemize}
\end{lemma}

It is natural to consider the following partial version of Bongartz (co-)completions.

\begin{definition}\label{def:parBgz}
Let $\xx\in\CCx_*(\seeds)$.
The \emph{partial Bongartz (co-)completions} are defined to be the partial cluster
\begin{equation}\label{eq:parBgz2}
    \Bgzpm{\xx}{\vv} \colon =\bigcap_{\xx\subset\ol{\xx}}\Bgzpm{\ol{\xx}}{\vv},
\end{equation}
where $\ol{\xx}$ runs through all clusters containing $\xx$.
\end{definition}

A consequence of \Cref{lem:Bong} is the following rank property for partial Bongartz (co-)completions.

\begin{corollary}\label{cor:rank}
For any partial clusters $\xx$ and $\vv$, $|\Bgzp{\xx}{\vv}|=|\Bgzm{\xx}{\vv}|$.
\end{corollary}

\begin{proof}
Let $\ol{\xx}$ be any cluster containing $\xx$.
\Cref{lem:Bong} provides the following identification
\[
\CEG^{\ol{\xx}}(\seeds,\vv)\cong\CEG^{\Bgzp{\ol{\xx}}{\vv}\backslash\vv}(\seeds_\vv),
\]
sending $\ol{\vv}\in\uCEG(\seeds,\vv)_0$ to $\ol{\vv}\backslash\vv\in\uCEG(\seeds_\vv)_0$.
Moreover, $\Bgzm{\ol{\xx}}{\vv}\backslash\vv$ is the unique sink of $\CEG^{\Bgzp{\ol{\xx}}{\vv}\backslash\vv}(\seeds_\vv)$.

Under \Cref{ass0}, the reduced cluster pattern $\seeds_\vv$ admits a green-to-red sequence so that the DT-transformation(cf. \Cref{def:DT}) is defined, denoted by $\DT{(-)_\vv}$.
It sends variables in $\Bgzp{\ol{\xx}}{\vv}\backslash\vv$ to $\Bgzm{\ol{\xx}}{\vv}\backslash\vv$ for any $\ol{\xx}$.

Restricting the bijection $\DT{(-)_\vv}$ to the partial Bongartz completion $\Bgzp{\xx}{\vv}\backslash\vv$, we have
\begin{gather}\label{eq:Bgz-DT}
    \Bgzp{\xx}{\vv}\backslash\vv
    \xlongleftrightarrow[1:1]{\DT{(-)_\vv}}
    {\Bgzm{{\xx}}{\vv}\backslash\vv}.
\end{gather}
Hence, $|\Bgzp{\xx}{\vv}\backslash\vv|=|\Bgzm{\xx}{\vv}\backslash\vv|$.
As $\vv\subset\Bgzpm{\xx}{\vv}$, the claim follows.
\end{proof}

\subsection{Parallel part of $\xevo$-evo flows}\label{sec:pevo}\

Let $\xevo=\sum_{i=1}^k\coex_i x_i\in\CCx(\seeds)$ be a point  with $\supp(\xevo)=\xx$ and $\ol{\xx}$ any cluster containing $\xevo$.
We extend $\flowdu{\xevo}$ to any partial cluster $\vv$ as follows.

\begin{definition}\label{def:parallel}
The \emph{parallel part of upward/downward $\xevo$-evo flow} is an operator
\begin{equation}\label{eq:flow}
\begin{array}{rcl}
\flowdu{\xevo}:
\CCx_*(\seeds) &\to    &\sCpx(\seeds)\\
\vv            &\mapsto           &\sum_{i=1}^k \coex_i\left[\normal_\pm(\tgvec^{\Bgzpm{\ol{\xx}}{\vv}}_{x_i})-\normal_\mp(\tgvec^{\Bgzpm{\ol{\xx}}{\vv}}_{x_i})\right].
\end{array}
\end{equation}
sending $\vv$ to a vector in $\RR\Bgzpm{\ol{\xx}}{\vv}$.
It is naturally continuous with respect to $\xevo$.
The continuity has many implications, e.g. an analogue of \Cref{lem:Bgz-mut}.
\end{definition}

\begin{remark}\label{rem:dege-evo}
Recall the fact from \cite[\S2.7]{Nak} that,
\begin{align}\label{eq:g=e}
[\widetilde{\gvec}^{\vv}_x]_- =0
 \iff x\in\vv,\qquad \widetilde{\gvec}^{\vv}_x=\bfit e_x\ \text{in that case};\\
\label{eq:g=-e}
[\widetilde{\gvec}^{\vv}_x]_+ =0
 \iff \DT{x}\in\vv,\qquad
 \widetilde{\gvec}^{\vv}_x=-\bfit e_{\DT{x}}\ \text{in that case}.
\end{align}
When the two special occasions above occur for some $\tgvec^{\Bgzpm{\ol{\xx}}{\vv}}_{x_i}$, 
$\flowdu{\xevo}$ is not necessarily a direction in $\dCpx(\seeds)$.
For instance, take $\xx=\supp(\xevo)$.
As $x\in\Bgzp{\ol\xx}{\xx}$ and $\DT{x}\in\Bgzp{\ol\xx}{\DT{\xx}}$ for any $x\in\supp\xevo$, we have
\begin{gather}\label{eq:ex-starX}
\begin{cases}
  \flowu{\xevo}(\xx)=\sum_{i=1}^k a_i \normal_+(\bfit{e}_{x_i})=\xevo,\\  
  \flowu{\xevo}(\DT{\xx})=-\sum_{i=1}^k a_i \normal_-(-\bfit{e}_{\DT{x_i}})=-\DT{\xevo}.
\end{cases}
\end{gather}
That is something we will remedy in \Cref{sec:evoCC}.
\end{remark}

We finish this subsection by showing the well-definedness of $\flowdu{\xevo}$.

\begin{lemma}\label{lem:well-define}
For any cluster variable $x$, the vectors
\begin{gather}\label{eq:Bgzpps}
    \normal_\pm(\tgvec^{\Bgzpm{\ol{\xx}}{\vv}}_x)
\end{gather}
do not depend on the choice of the cluster $\ol{\xx}$ containing $x$.
\end{lemma}

\begin{proof}
Let $\ol{\xx},\ol{\xx}'$ be two clusters containing $x$.
The subgraph $\uCEG(\seeds,x)$ of clusters containing $x$ is connected since it is the exchange graph of the reduced cluster pattern $\seeds_x$.
Hence it suffices to consider adjacent clusters $\ol{\xx}'=\mu_{x_k}(\ol{\xx})$ with $x_k\ne x$.
We only need to show that
\begin{gather}\label{eq:tech1}
    \tgvec^{\Bgzpm{\ol{\xx}}{\vv}}_x=\tgvec^{\Bgzpm{\ol{\xx}'}{\vv}}_x
\end{gather}
for any $x\in\ol{\xx}\cap\ol{\xx}'$.
Denote by $\ol{\vv}=\Bgzpm{\ol{\xx}}{\vv}$ and $\ol{\vv}'=\Bgzpm{\ol{\xx}'}{\vv}$.
By \Cref{lem:Bgz-mut}, we have
\[
\ol{\vv}'=\ol{\vv}~\text{or}~\mu_{v_0}(\ol{\vv})~\text{for some}~v_0\in \ol{\vv} \backslash\vv.
\]
When $\ol{\vv}=\ol{\vv}'$, \eqref{eq:tech1} follows immediately.

When $\ol{\vv}'=\mu_{v_0}(\ol{\vv})$ for some $v_0\in \ol{\vv} \backslash\vv$, 
denote by 
\[
B_{\ol{\vv}}=(b_{v_i, v_j})_{v_i,v_j\in\ol{\vv}},\quad
\Cmat_{\ol{\vv}}^{\ol{\xx}}=(c_{x_i,v_j})_{x_i\in\ol{\xx},v_j\in\ol{\vv}}\quad\text{and}\quad
\Cmat_{\ol{\vv}'}^{\ol{\xx}}=(c'_{x_i,v_j })_{x_i\in\ol{\xx},v_j\in\ol{\vv}'}
\]
By \eqref{eq:mutC}, the mutation at $v_0\in\ol{\vv}$ gives
\begin{gather}\label{eq:C0}
     \Cmat_{\ol{\vv}'}^{\ol{\xx}} =
    \Cmat_{\ol{\vv}}^{\ol{\xx}}J_{v_0} +
    \Cmat_{\ol{\vv}}^{\ol{\xx}} E_{v_0 v_0} [B_{\ol{\vv}}]_+ +
    [-\Cmat_{\ol{\vv}}^{\ol{\xx}}]_+ E_{v_0 v_0} B_{\ol{\vv}}.
\end{gather}


For any $x\in\ol{\xx}\cap\ol{\xx}'$, the relation \eqref{eq:C0} restricts to the $(x,v_0)$-entry gives
\begin{gather}\label{eq:c0}
    c'_{x,v_0} =
    -c_{x, v_0} + c_{x,v_0} [b_{v_0, v_0}]_+  + [-c_{x,v_0}]_+ b_{v_0, v_0}
    =-c_{x, v_0}
\end{gather}
By the duality \eqref{eq:dual3}, we have
\[
c'_{x,v_0}=\ubf{e}_{x}\cvec_{v_0,\ol{\vv}'}^{\ol{\xx}}
=\ubf{e}_{v_0}\widetilde{\gvec}_{x}^{\ol{\vv}'}
=\ubf{e}_{x}\cvec_{v_0,\ol{\vv}'}^{\ol{\xx}'}.
\]
Hence $c'_{x,v_0}$ is also an entry of $\cvec_{v_0,\ol{\vv}'}^{\ol{\xx}'}$.
But $\cvec_{v_0,\ol{\vv}'}^{\ol{\xx}'}$ and $\cvec_{v_0,\ol{\vv}}^{\ol{\xx}}$ share the same sign.
Combining with \eqref{eq:c0}, we have 
\begin{gather}\label{eq:c1}
    c_{x,v_0}=0=c'_{x,v_0}.   
\end{gather}
Restricting \eqref{eq:C0} to the $(x,v)$-entry for any $v\ne v_0$, we have
\begin{gather}\label{eq:c2}
     c'_{x,v} =
    c_{x, v} + c_{x,v_0} [b_{v_0, v}]_+ + [-c_{x,v_0}]_+ b_{v_0, v}
    =c_{x,v}.
\end{gather}
Consequently we deduce that
\[
\tgvec^{\ol\vv}_{x}\xlongequal{\eqref{eq:dual3}}
(\ubf{e}_{x}\Cmat^{\ol{\xx}}_{\ol{\vv}})^T\xlongequal{\eqref{eq:c1}+\eqref{eq:c2}}
(\ubf{e}_{x}\Cmat^{\ol{\xx}}_{\ol{\vv}'})^T\xlongequal{\eqref{eq:dual3}}
\tgvec^{\ol{\vv}'}_{x}.
\]
\eqref{eq:tech1} holds as required.
\end{proof}

As a consequence of \Cref{lem:well-define}, we see that $\flowdu{\xevo}(\vv)$ in fact lives in
\[
    \RR\Bgzpm{\xx}{\vv}=\bigcap_{\xx\subset\ol{\xx}}\RR\Bgzpm{\ol{\xx}}{\vv}.
\]
In other words, we have
\begin{gather}\label{eq:suppevo}
    \supp(\flowdu{\xevo}(\vv))\subset\Bgzpm{\xx}{\vv}.
\end{gather}
We proceed to show that $\Bgzpm{\xx}{\vv}$ is precisely the support for $\flowdu{\xevo}(\vv)$ at $\vv$
in a suitable sense.

\subsection{$\xevo$-Supports and the rank invariant}\label{sec:Xsupp}\


\begin{definition}\label{def:Xsupp}
The \emph{upward/downward $\xevo$-support} is an operation
\begin{equation}\label{eq:Xsupp-V}
\begin{array}{rcl}
\suppdu{\xevo}:\CCx_*(\seeds) &\to    &\CCx_*(\seeds)\\
\vv                           &\mapsto&\vv\cup\supp(\flowdu{\xevo}(\vv)),
\end{array}
\end{equation}
\end{definition}

Note that $\supp(\flowdu{\xevo}(\vv))$ only depends on $\xx=\supp(\xevo)$ so we may use short version $\xx$ instead of $\xevo$. 
\begin{gather}\label{eq:simplified}
    \suppdu{\xevo}(\vv)=\suppdu{\xx}(\vv)
\end{gather}
The $\xevo$-support is bounded below by both $|\vv|$ and $|\xx|$.
The latter lower bound is not clear from definition.
\begin{lemma}\label{lem:bbelow-xx}
    For any $\xx,\vv\in\CCx_*(\seeds)$, we have
    \begin{gather}\label{eq:bbelow-xx}
        |\suppdu{\xx}(\vv)|\ge|\xx|.
    \end{gather}
\end{lemma}

\begin{proof}
We claim that
\[
\supp(\flowdu{\xevo}(\vv))
=\bigcup_{x\in\xx}\supp\left(\tgvec^{\ol\vv}_x\right),
\qquad\ol{\vv}=\Bgzpm{\ol{\xx}}{\vv},
\quad \xx\subset\ol{\xx}.
\]
Indeed, by definition we have the inclusion
\[
\supp(\flowdu{\xevo}(\vv))
\subset\bigcup_{x\in\xx}\supp\left(\tgvec^{\ol\vv}_x\right).
\]
Conversely,
suppose $v\in\supp\left(\tgvec^{\ol\vv}_x\right)$ for some $x\in\xx$,
The row sign-coherence \eqref{eq:g-sc} ensures that all nonzero contributions from $\{\tgvec^{\ol\vv}_x\}_{x\in\xx}$ to the $v$-coordinate in $\flowdu{\xevo}(\vv)$ have the same sign; hence no cancellation occurs and
\begin{gather}\label{eq:supp-g}
    \supp\left(\tgvec^{\ol\vv}_x\right)\subset\supp(\flowdu{\xevo}(\vv)),\qquad\forall x\in\xx.
\end{gather}
The claim holds as desired.
Since $\{\tgvec^{\ol\vv}_x\}_{x\in\xx}$ are linearly independent, we have
\[
|\xx|\le
|\bigcup_{x\in\xx}  \supp\left(\tgvec^{\ol\vv}_x\right)|=
|\supp(\flowdu{\xevo}(\vv))|\le
|\suppdu{\xx}(\vv)|.\qedhere
\]
\end{proof}

The partial Bongartz (co-)completions and $\xevo$-supports coincide when $\xx=\supp(\xevo)$ is a cluster.
This is indeed the case in general.
We first prove a technical lemma.

\begin{lemma}\label{lem:star-correspondence}
Let $\xx,\vv\in\CCx_*(\seeds)$ such that $\vv=\suppdu{\xx}(\vv)$.
Then the map
\[
\begin{array}{rcl}
\Phi:\CEG_0(\seeds,\xx)&\longrightarrow&\CEG_0(\seeds,\vv)\\
\ol{\xx}&\longmapsto&\Bgzpm{\ol{\xx}}{\vv}
\end{array}
\]
is surjective with a section
\[
\begin{array}{rcl}
\Psi:\CEG_0(\seeds,\vv)&\longrightarrow&\CEG_0(\seeds,\xx)\\
\ol{\vv}&\longmapsto&
\tBgzpm{\ol{\vv}}{\xx},
\end{array}
\]
where $\widetilde{\operatorname{Bgz}}^\updownarrow$ is the transposed Bongartz (co-)completion, cf. \Cref{def:bong}.
Moreover, if $|\suppdu{\xx}(\vv)|=|\xx|$, then $\Phi$ and $\Psi$ give the 1-1 correspondence between $\CEG_0(\seeds,\vv)$ and $\CEG_0(\seeds,\xx)$.
\end{lemma}

\begin{proof}
We treat the upward and downward cases simultaneously,
writing $\varepsilon=1$ for $\uparrow$ and
$\varepsilon=-1$ for $\downarrow$.

By \Cref{def:bong}, the map $\Phi$ is well-defined.
We claim that
\begin{gather}\label{eq:star-vanishing}
    \supp(\tgvec^{\ol{\vv}}_x)\subset\vv,
    \qquad
    \forall x\in\xx,\quad
    \forall\ol{\vv}\supset\vv.
\end{gather}
Indeed, taking a completion $\xx\subset\ol{\xx}$ this holds for $\ol{\vv}_0=\Bgzpm{\ol{\xx}}{\vv}$ by \eqref{eq:supp-g}.
By \cite[Thm.~10]{CL2}, every cluster containing $\vv$ is connected to $\ol{\vv}_0$ by a sequence of mutations outside $\vv$.

Suppose $\supp(\tgvec^{\ol{\vv}}_x)\subset\vv$ and $\ol{\vv}'=\mu_v(\ol{\vv})$ for $v\in\ol{\vv}\backslash\vv$.
By definition \eqref{eq:dual3},
the $v$-entry of $\ubf{e}_x\Cmat^{\ol{\xx}}_{\ol{\vv}}$ vanishes.
The mutation formula \eqref{eq:mutC} then shows that this row remains unchanged under the mutation at $v$.
Thus \eqref{eq:star-vanishing} follows by induction.

Now take any $\ol{\vv}\in\CEG_0(\seeds,\vv)$ and put
\[
    \ol{\xx}'=
    \tBgzpm{\ol{\vv}}{\xx}.
\]
Its existence and uniqueness follow from \Cref{lem:Bong} and \Cref{ass0}, applied to the cluster pattern with initial seed $(\ol{\vv},B_{\ol{\vv}}^T)$.
By \eqref{eq:star-vanishing},
\begin{gather}\label{eq:star-zero}
    (\widetilde{\Gmat}^{\ol{\vv}}_{\ol{\xx}'})_{v,x}=0,
    \qquad
    \forall v\in\ol{\vv}\backslash\vv,\quad
    \forall x\in\xx.
\end{gather}

By the defining property of Bongartz
(co-)completion with respect to transposed $C$-matrices,
\[
    \varepsilon\,
    \widetilde{\cvec}^{\ol{\vv}}_{y,\ol{\xx}'}
    \in\ZZp^n,
    \qquad
    \forall y\in\ol{\xx}'\backslash\xx.
\]
Using tropical duality \eqref{eq:dual2} and the positivity of the skew-symmetrizer, we obtain
\begin{gather}\label{eq:star-inverse}
    \varepsilon(\widetilde{\Gmat}^{\ol{\vv}}_{\ol{\xx}'})^{-1}_{y,v}\ge0,
    \qquad
    \forall y\in\ol{\xx}'\backslash\xx,\quad
    \forall v\in\ol{\vv}.
\end{gather}
Fix $v\in\ol{\vv}\backslash\vv$.
By row sign-coherence \eqref{eq:g-sc}, the $v$-row of $\widetilde{\Gmat}^{\ol{\vv}}_{\ol{\xx}'}$ is sign-coherent.
If we have $\varepsilon (\widetilde{\Gmat}^{\ol{\vv}}_{\ol{\xx}'})_{v,y}\le0$ for every $y\in\ol{\xx}'$, then
\eqref{eq:star-zero} and \eqref{eq:star-inverse} imply
\[
    1=\left(\widetilde{\Gmat}^{\ol{\vv}}_{\ol{\xx}'}(\widetilde{\Gmat}^{\ol{\vv}}_{\ol{\xx}'})^{-1}\right)_{v,v}
     =\sum_{y\in\ol{\xx}'\backslash\xx}
       (\widetilde{\Gmat}^{\ol{\vv}}_{\ol{\xx}'})_{v,y}(\widetilde{\Gmat}^{\ol{\vv}}_{\ol{\xx}'})^{-1}_{y,v}
     \le0,
\]
a contradiction. Therefore
\[
    \varepsilon (\widetilde{\Gmat}^{\ol{\vv}}_{\ol{\xx}'})_{v,y}\ge0,
    \qquad
    \forall v\in\ol{\vv}\backslash\vv,\quad
    \forall y\in\ol{\xx}'.
\]
By definition \eqref{eq:dual3},
\[
    \Cmat^{\ol{\xx}'}_{\ol{\vv}}=(\widetilde{\Gmat}^{\ol{\vv}}_{\ol{\xx}'})^T.
\]
Thus
\[
    \varepsilon\,
    \cvec^{\ol{\xx}'}_{v,\ol{\vv}}
    \in\ZZp^n,
    \qquad
    \forall v\in\ol{\vv}\backslash\vv.
\]
Since $\vv\subset\ol{\vv}$, the uniqueness
of Bongartz (co-)completion gives
\begin{gather}\label{eq:star-bongartz-u}
    \ol{\vv}=\Bgzpm{\ol{\xx}'}{\vv}.
\end{gather}
Consequently,
\[
    \Phi\circ\Psi
    =\id_{\CEG_0(\seeds,\vv)},
\]
which proves the surjectivity of $\Phi$.
Conversely, suppose $|\vv|=|\suppdu{\xx}(\vv)|=|\xx|$.
For any $\ol{\xx}\supset\xx$, put
$\ol{\vv}=\Bgzpm{\ol{\xx}}{\vv}$ and
$\widetilde{\Gmat}^{\ol{\vv}}_{\ol{\xx}}$.
By \eqref{eq:star-vanishing}, after ordering
the variables with $\vv$ and $\xx$ first,
\[
    \widetilde{\Gmat}^{\ol{\vv}}_{\ol{\xx}}=
    \begin{pmatrix}
        A&B\\
        0&M
    \end{pmatrix},
\]
where $A$ and $M$ are square and invertible
by \eqref{eq:det1}.

Since $\ol{\vv}=\Bgzpm{\ol{\xx}}{\vv}$, the defining sign condition and \eqref{eq:dual3} give
\[
    \varepsilon (M)_{v,x}\ge0
\]
entrywise.
Similar to the previous argument, 
\eqref{eq:dual2} and column sign-coherence \eqref{eq:c-sc} show that every row of $M^{-1}$ is sign-coherent.
A row of $M^{-1}$ with the opposite sign would make the corresponding diagonal entry of $I=M^{-1}M$ non-positive.
This forces the following entrywise positivity for $M^{-1}$.
\[
    \varepsilon (M^{-1})_{v,x}\ge0.
\]
Using the block form of $(\widetilde{\Gmat}^{\ol{\vv}}_{\ol{\xx}})^{-1}$ and \eqref{eq:dual2} again, we conclude that
\[
    \varepsilon\,
    \widetilde{\cvec}^{\ol{\vv}}_{y,\ol{\xx}}
    \in\ZZp^n,
    \qquad
    \forall y\in\ol{\xx}\backslash\xx.
\]
Thus $\ol{\xx}$ is the dual Bongartz
(co-)completion of $\xx$ with respect to
$\ol{\vv}$, i.e.,
\[
    \Psi\circ\Phi=\id_{\CEG_0(\seeds,\xx)}.
\]
Hence the correspondence is bijective.
\end{proof}

\begin{theorem}\label{thm:X=xx}
For any $\xx,\vv\in\CCx_*(\seeds)$, we have
\begin{gather}\label{eq:X=xx}
    \suppdu{\xx}(\vv)=\Bgzpm{\xx}{\vv}
\end{gather}
Moreover, $|\suppu{\xx}(\vv)|=|\suppd{\xx}(\vv)|\ge|\xx|$.
\end{theorem}

\begin{proof}
Without loss of generality, we may assume
\[
    \vv=\suppdu{\xx}(\vv).
\]
By \Cref{lem:star-correspondence}, the map
\[
    \ol{\xx}\longmapsto
    \Bgzpm{\ol{\xx}}{\vv}
\]
is surjective from $\CEG_0(\seeds,\xx)$ onto
$\CEG_0(\seeds,\vv)$.
Consequently, by \Cref{def:parBgz},
\begin{gather*}
    \Bgzpm{\xx}{\vv}=
    \bigcap_{\xx\subset\ol{\xx}}\Bgzpm{\ol{\xx}}{\vv}=
    \bigcap_{\vv\subset\ol{\vv}}\ol{\vv}=
    \vv=\suppdu{\xx}(\vv).
\end{gather*}
Here the third equality follows because every variable outside $\vv$ can be removed by a mutation fixing $\vv$.
This proves \eqref{eq:X=xx}.
The equality
\[
    |\suppu{\xx}(\vv)|
    =
    |\suppd{\xx}(\vv)|
\]
now follows from \Cref{cor:rank}.
\end{proof}

\begin{remark}\label{rem:union}
While $\Bgzpm{\xx}{\vv}$ is defined in an intrinsic way,
the definition of $\suppdu{\xx}(\vv)$ has its own advantage.
For instance, row sign-coherence \eqref{eq:g-sc} ensures that no coordinate cancellation occurs. 
Hence one has
\[
    \suppdu{\xx}(\vv)=\bigcup_{x\in\xx} \suppdu{x}(\vv)
\]
or $\Bgzpm{\xx}{\vv}=\bigcup_{x\in\xx}\Bgzpm{x}{\vv}$, which is not obvious from \eqref{eq:parBgz2}.
\end{remark}

\section{\texorpdfstring{$\xevo$-\EVO flows and $\xevo$-foliations}{X-evo flows and X-foliations}}\label{sec:evoCC}

In this section, we modify $\flowdu{\xevo}$ to be a direction field, called $\xevo$-\EVO flow,
so that it induces a nice (oriented) foliation on the cluster complex.

\subsection{A modification of $\xevo$-\EVO flow}\label{sec:evo}\

Let $\xevo=\sum_{i=1}^k \coex_i x_i$ and $\ccpt{V}$ be points in $\CCx(\seeds)$.
Set $\xx=\supp(\xevo)$ and $\vv=\supp(\ccpt{V})$.
Let $\ol{\xx}$ be a cluster containing $\xx$ .

\begin{definition}\label{def:Evo}
The \emph{upward/downward $\xevo$-\EVO flow} is a one-sided piecewise-defined direction field, i.e. a function
\begin{equation}\label{eq:UDflow}
\begin{array}{rcl}
\flowDU{\xevo}:\CCx(\seeds) &\to     &\dCpx(\seeds)  \\
\ccpt{V}                    &\mapsto &\flowdu{\xevo}(\vv)-\csum(\flowdu{\xevo}(\vv))\cdot\ccpt{V}.
\end{array}
\end{equation}
sending a point $\ccpt{V}$ to a direction in $\dcc(\Bgzpm{{\xx}}{\vv})$.
\end{definition}

\begin{lemma}\label{lem:sing}
$\flowDU{\xevo}(\ccpt{V})$ is degenerated (i.e. $\flowDU{\xevo}(\ccpt{V})=\bfit{0}$) if and only if $\ccpt{V}=\xevo$ or $\Xevo$.
\end{lemma}

\begin{proof}
We prove the upward statement and the downward version is similar.

Applying \Cref{def:Evo} to \eqref{eq:ex-starX} shows that  
\begin{gather}
    \flowU{\xevo}(\xevo)=\bfit{0}=\flowU{\xevo}(\Xevo).
\end{gather}
By row sign-coherence \eqref{eq:g-sc}, all nonzero contributions to each coordinate of $\flowu{\xevo}(\vv)$ have the same sign, which prevents cancellation.

Every selected $g$-vector is nonzero, and hence
\[
\flowu{\xevo}(\vv)\ne\bfit0.
\]
Suppose $\flowU{\xevo}(\ccpt{V})=\bfit{0}$. 
We have
\begin{gather}
    \flowu{\xevo}(\vv) =\csum(\flowu{\xevo}(\vv))\ccpt{V}, 
    \quad\csum(\flowu{\xevo}(\vv))\ne0.
\end{gather}
If $\csum(\flowu{\xevo}(\vv))>0$, then
\[
[\flowu{\xevo}(\vv)]_-=\bfit0.
\]
Since cancellation is impossible and all $a_i$ are positive, we obtain
\begin{gather}
   [\tgvec^{\Bgzp{\ol\xx}{\vv}}_{x_i}]_-=\bfit{0},\quad\forall~i.
\end{gather}
\Cref{eq:g=e} therefore gives
\begin{gather}
    \tgvec^{\Bgzp{\ol\xx}{\vv}}_{x_i}=\bfit{e}_{x_i},\quad
    \flowu{\xevo}(\vv)=\xevo,\quad
    \csum(\flowu{\xevo}(\vv))=1.
\end{gather}
It follows that $\ccpt{V}=\xevo$.
If $\csum(\flowu{\xevo}(\vv))<0$, the same argument instead gives
\begin{gather}
   [\tgvec^{\Bgzp{\ol\xx}{\vv}}_{x_i}]_+=\bfit{0},\quad\forall~i.
\end{gather}
By \eqref{eq:g=-e}, 
\begin{gather}
    \tgvec^{\Bgzp{\ol\xx}{\vv}}_{x_i}=-\bfit{e}_{\DT{x}_i},\quad
    \flowu{\xevo}(\vv)=-\Xevo,\quad
    \csum(\flowu{\xevo}(\vv))=-1.
\end{gather}
Thus $\ccpt{V}=\Xevo$. 
\end{proof}

\begin{remark}\label{rem:Xsupp}
When we consider the $\xevo$-evo flow, the $\xevo$-support (cf. \Cref{def:Xsupp}) for $\supp(\ccpt{V})$ in fact illustrates the support of the nearby perturbed points of $\ccpt{V}$ along $\flowDU{\xevo}$, i.e.,
\begin{gather}\label{eq:suppEvo}
\suppdu{\xevo}(\supp(\ccpt{V}))=\supp(\ccpt{V}+\varepsilon\flowDU{\xevo}(\ccpt{V})),
\qquad     0<\varepsilon\ll1.
\end{gather}
\end{remark}

When $\vv$ is a cluster, we have $\suppu{\xevo}(\vv)=\vv=\suppd{\xevo}(\vv)$.
In such a case, $\flowu{\xevo}(\vv)=-\flowd{\xevo}(\vv)$ so $\flowU{\xevo}(\ccpt{V})=-\flowD{\xevo}(\ccpt{V})$ are compatible.
In general, the compatibility of $\flowDU{\xevo}(\vv)$ depends on $\suppdu{\xevo}(\vv)$ as follows.

\begin{lemma}\label{lem:c-coincide}
For any $\ccpt{V}$, we have
\begin{gather}\label{eq:c-coincide}
    \suppu{\xevo}(\vv)=\suppd{\xevo}(\vv)\iff
    \vv=\suppdu{\xevo}(\vv)\iff
    \flowU{\xevo}(\ccpt{V})=-\flowD{\xevo}(\ccpt{V}).
\end{gather}
\end{lemma}

\begin{proof}
First suppose that 
\[
\suppu{\xevo}(\vv)=\suppd{\xevo}(\vv).
\]
Both $\Bgzp{\ol{\xx}}{\vv}$ and $\Bgzm{\ol{\xx}}{\vv}$ contain this common partial cluster.
Connectivity of its reduced exchange graph provides a mutation path
\begin{gather}\label{eq:Bgzp-m}
    \Bgzp{\ol{\xx}}{{\vv}}=
    \ww_0\xrightarrow{\mu_{w_0}}\ldots\xrightarrow{\mu_{w_{N-1}}}\ww_N=\Bgzm{\ol{\xx}}{\vv}.
\end{gather}
fixing every variable in $\suppdu{\xevo}(\vv)$.
By \Cref{def:Xsupp} and no-cancellation property, the initial columns satisfy
\begin{gather}
    \ubf{e}_{w}\tgvec^{\ww_0}_{x_i}=
    \ubf{e}_{x_i}\cvec^{\ol{\xx}}_{w,\ww_0}=0~\text{for every}~i~\text{and}~w\in\ww_0\backslash\suppu{\xevo}(\vv).
\end{gather}
The vanishing holds at \emph{all} positions outside $\suppu{\xevo}(\vv)$.
Moreover, if it holds at $\ww_j$, then the entry at the mutated position $w_j$ is zero.
The mutation formula \eqref{eq:mutC} consequently gives
\begin{gather}
    \ubf{e}_{x_i}\Cmat^{\ol{\xx}}_{\ww_{j+1}}=\ubf{e}_{x_i}\Cmat^{\ol{\xx}}_{\ww_j}.
\end{gather}
Indeed, every correction term contains the zero entry at $w_j$, and changing the sign of that entry leaves it zero. 
Thus the same vanishing condition holds at $\ww_{j+1}$.
Induction proves equality of the selected rows throughout the path \eqref{eq:Bgzp-m}.

For any $v\in\suppu{\xevo}(\vv)\backslash\vv$, the definition of Bongartz (co-)completions gives
\[
    \cvec^{\ol{\xx}}_{v,\Bgzp{\ol{\xx}}{\vv}}\in\ZZ^n_{\ge0},\quad
    \cvec^{\ol{\xx}}_{v,\Bgzm{\ol{\xx}}{\vv}}\in\ZZ^n_{\le0}.
\]
The row equality of $C$-matrices just proved imply, for every $i$,
\begin{gather}
    0\le \ubf{e}_{x_i}\cvec^{\ol{\xx}}_{v,\Bgzp{\ol{\xx}}{\vv}}
    =\ubf{e}_{x_i}\cvec^{\ol{\xx}}_{v,\Bgzm{\ol{\xx}}{\vv}}\le0.
\end{gather}
All selected entries therefore vanish.
The $v$-coordinate of $\flowu{\xevo}(\vv)$ is zero, contrary to $v\in\suppu{\xevo}(\vv)\backslash\vv$.
We conclude that
\begin{gather}
    \suppu{\xevo}(\vv)=\vv=\suppd{\xevo}(\vv).
\end{gather}
Next suppose that
\begin{gather}
    \vv=\suppu{\xevo}(\vv)
    \quad\text{or}\quad
    \vv=\suppd{\xevo}(\vv).
\end{gather}
So in the corresponding $g$-vectors $\tgvec^{\Bgzp{\ol\xx}{\vv}}_{x_i}$ or $\tgvec^{\Bgzm{\ol\xx}{\vv}}_{x_i}$ all rows vanish outside $\vv$.
Choose a mutation path between $\Bgzp{\ol{\xx}}{\vv}$ and $\Bgzm{\ol{\xx}}{\vv}$ fixing $\vv$ and apply the preceding induction again, 
starting at whichever endpoint has this vanishing property.
It follows that
\begin{gather}
    \tgvec^{\Bgzp{\ol{\xx}}{\vv}}_{x_i}
    =(\ubf{e}_{x_i}\Cmat^{\ol{\xx}}_{\Bgzp{\ol{\xx}}{\vv}})^T
    =(\ubf{e}_{x_i}\Cmat^{\ol{\xx}}_{\Bgzm{\ol{\xx}}{\vv}})^T
    =\tgvec^{\Bgzm{\ol{\xx}}{\vv}}_{x_i}.
\end{gather}
So by \Cref{def:parallel}, we have
\[
    \flowu{\xevo}(\vv)=-\flowd{\xevo}(\vv).
\]
Both supports equal to $\vv$, and subtracting the coefficient sums times $\ccpt{V}$ gives
\begin{gather}
    \flowU{\xevo}(\vv)=-\flowD{\xevo}(\vv).
\end{gather}
Finally, suppose $\flowU{\xevo}(\ccpt{V})=-\flowD{\xevo}(\ccpt{V})$. 
Then their coordinates outside $\vv$ are both non-negative and opposite to each other, hence zero.
\Cref{eq:suppEvo} gives both supports equal to $\vv$. 
\end{proof}

If both upward and downward $\xevo$-evo flows supports in $\vv$, then the direction field on the relative interior
decompose it into straight lines:
\begin{itemize}
    \item When $\csum(\flowu{\xevo}(\vv))=0$, they are parallel lines;
    \item otherwise they lie on lines through $\flowu{\xevo}(\vv)/\csum(\flowu{\xevo}(\vv))$.
\end{itemize}

\begin{definition}
$\flowDU{\xevo}$ induce a \emph{local linear $\xevo$-foliation} on each $\cc(\vv)^\circ$ with $\suppu{\xevo}(\vv)=\suppd{\xevo}(\vv)$, i.e., there is a decomposition
\begin{gather}\label{eq:local-foli}
\cc(\vv)^\circ\backslash\{\xevo,\Xevo\}=\bigsqcup_{l~\text{local $\xevo$-leaves}} l.
\end{gather}
\end{definition}

\begin{remark}\label{rem:dege-foli}\
\begin{enumerate}
    \item If $\suppu{\xevo}(\vv)=\vv=\suppd{\xevo}(\vv)$, then the transport observation in the proof of \Cref{lem:c-coincide} shows that the selected rows are identical in every cluster containing $\vv$.
    Hence for $\vv\subset\vv'$ the parallel part of fields in $\vv'$ agree with that at $\vv$, and the radial modifications restrict compatibly to the face $\cc(\vv)^\circ$.
    
    \item If $\suppu{\xevo}(\vv)\ne\suppd{\xevo}(\vv)$, then both one-sided fields leave the relative interior of $\vv$.
\end{enumerate}
\end{remark}


\subsection{$\xevo$-Foliation and cell-crossing}\label{sec:Xfoli}\

Our goal is to glue the local $\xevo$-foliations on each cell and obtain a global foliation structure on $\CCx(\seeds)$ with two singularities.

\begin{lemma}\label{lem:c-crossing}
    Every regular point $\ccpt{V}\ne\xevo,\Xevo$ admits a unique outgoing and a unique incoming oriented local $\xevo$-leaf germ with respect to $\flowU{\xevo}$.
    Their supporting cells are respectively $\cc(\suppu{\xevo}(\vv))$ and $\cc(\suppd{\xevo}(\vv))$.

    Consequently, any regular local $\xevo$-leaf has a unique oriented continuation at each regular endpoint. 
    Singular endpoints terminate the regular leaf.
\end{lemma}

\begin{proof}
    Fix a cluster $\ol{\xx}$ containing $\supp\xevo$ and a regular point $\ccpt{V}_0$.
    Put $\vv_0=\supp(\ccpt{V}_0)$.
    We first construct the outgoing local $\xevo$-leaf germ.

    By \Cref{def:Xsupp},
    \begin{gather}
        \vv_0\subset\suppu{\xevo}(\vv_0)\subset\Bgzp{\ol{\xx}}{\vv_0}.
    \end{gather}
    Every $c$-vector of $\Bgzp{\ol{\xx}}{\vv_0}$ outside $\vv_0$ is non-negative hence so is every $c$-vector outside $\suppu{\xevo}(\vv_0)$.
    The Uniqueness of Bongartz completion gives
    \begin{gather}
        \Bgzp{\ol{\xx}}{\vv_0}=\Bgzp{\ol{\xx}}{\suppu{\xevo}(\vv_0)}.
    \end{gather}
    Applying this to \Cref{def:parallel} and \Cref{def:Xsupp}, we obtain
    \begin{gather}\label{eq:evo1}
        \flowu{\xevo}(\suppu{\xevo}(\vv_0))=\flowu{\xevo}(\vv_0),\\
        \suppu{\xevo}(\suppu{\xevo}(\vv_0))=\suppu{\xevo}(\vv_0).
    \end{gather}
    \Cref{lem:c-coincide} shows that 
    \begin{gather}
        \suppd{\xevo}(\suppu{\xevo}(\vv_0))=\suppu{\xevo}(\suppu{\xevo}(\vv_0)).
    \end{gather}
    Thus its relative interior carries a local $\xevo$-foliation. 
    For $0<\varepsilon\ll1$, \eqref{eq:suppEvo} gives
    \begin{gather}
        \ccpt{V}_0+\varepsilon\flowU{\xevo}(\ccpt{V}_0)\in\cc(\suppu{\xevo}(\vv_0))^\circ.
    \end{gather}
    \Cref{eq:evo1} above and \Cref{def:Evo} provide a more precise identity
    \begin{align*}
        &\flowU{\xevo}(\ccpt{V}_0+\varepsilon\flowU{\xevo}(\ccpt{V}_0))\\
        =&\flowu{\xevo}(    \suppu{\xevo}(\vv_0)    )   -  
        \csum(  
        \flowu{\xevo}(\suppu{\xevo}(\vv_0))
        )\cdot(    \ccpt{V}_0+\varepsilon\flowU{\xevo}(\ccpt{V}_0) 
        )\\
        =&\left[    \flowu{\xevo}(\vv_0)    -   \csum(\flowu{\xevo}(\vv_0))\ccpt{V}_0    \right]    -   \varepsilon\csum(\flowu{\xevo}(\vv_0))  \flowU{\xevo}(\ccpt{V}_0)\\
        =&\left[
        1-\varepsilon\csum(\flowu{\xevo}(\vv_0)) 
        \right]
        \flowU{\xevo}(\ccpt{V}_0).
    \end{align*}
    The scalar factor is positive for sufficiently small $\varepsilon$.
    Since $\ccpt{V}_0$ is regular, the direction is nonzero.
    This constructs the outgoing straight local leaf near $\ccpt{V}_0$ and proves uniqueness within $\cc(\suppu{\xevo}(\vv_0))^\circ$.

    We proceed to show the uniqueness in different cells.
    Suppose that $\vv_0\subsetneqq\vv$ and a regular local $\xevo$-leaf $l$ in $\cc(\vv)^\circ$ issues from $\ccpt{V}_0$.
    Choose $\ccpt{V}_1\in l$ sufficiently near $\ccpt{V}_0$ and let
    \begin{gather}
        \ccpt{V}_t=(1-t)\ccpt{V}_0+t\ccpt{V}_1,\quad 0<t\le1.
    \end{gather}
    On this segment we have
    \begin{gather}
        \flowU{\xevo}(\ccpt{V}_1)=
        \flowu{\xevo}(\vv)-\csum(\flowu{\xevo}(\vv))\ccpt{V}_1=
        \lambda_1(\ccpt{V}_1-\ccpt{V}_0),
        \quad\lambda_1>0.
    \end{gather}
    By computation,
    \begin{gather}
        \flowU{\xevo}(\ccpt{V}_t)=
        \left[
        \lambda_1+(1-t)\csum(\flowu{\xevo}(\vv))
        \right]
        (\ccpt{V}_1-\ccpt{V}_0)
        =\colon\lambda_t(\ccpt{V}_1-\ccpt{V}_0),
    \end{gather}
    where $\lambda_t>0$.
    Taking $t\to 0^+$, we have
    \begin{gather}\label{eq:EvoV0}
        \flowu{\xevo}(\vv)-\csum(\flowu{\xevo}(\vv))\ccpt{V}_0=
        \lambda_0(\ccpt{V}_1-\ccpt{V}_0),
    \end{gather}
    where $\lambda_0\ge0$.
    An argument similar to \Cref{lem:sing} shows that $\lambda_0\ne0$. 
    Hence $\lambda_0>0$.
    
    For every $v\in\vv\setminus\vv_0$, the $v$-coordinate of
    $\ccpt{V}_0$ is zero, whereas that of $\ccpt{V}_1$ is positive.
    \eqref{eq:EvoV0} therefore shows that $v$-coordinate of $\flowu{\xevo}(\vv)$ is positive.
    By \Cref{def:parallel}, there exists $x_i\in\supp(\xevo)$ such that
    \[
    \ubf{e}_{v}\tgvec^{\Bgzp{\ol{\xx}}{\vv}}_{x_i}=\ubf{e}_{x_i}\cvec^{\ol{\xx}}_{v,\Bgzp{\ol{\xx}}{\vv}}>0.
    \]
    Sign-coherence implies
    \begin{gather}
        \cvec^{\ol{\xx}}_{v,\Bgzp{\ol{\xx}}{\vv}}\in\ZZp^n,\qquad
        v\in\vv\setminus\vv_0.
    \end{gather}
    For variables outside $\vv$, the same nonnegativity follows from the definition of $\Bgzp{\ol{\xx}}{\vv}$.
    Thus every $c$-vector of this cluster outside $\vv_0$ is nonnegative.
    Uniqueness of Bongartz completion now gives
    \[
        \Bgzp{\ol{\xx}}{\vv}
        =\Bgzp{\ol{\xx}}{\vv_0},
        \qquad
        \flowu{\xevo}(\vv)=\flowu{\xevo}(\vv_0).
    \]
    Because $\suppu{\xevo}(\vv)=\vv$, the parallel vector $\flowu{\xevo}(\vv)$ is supported on $\vv$.
    Consequently,
    \[
        \suppu{\xevo}(\vv_0)
        =
        \vv_0\cup\supp(\flowu{\xevo}(\vv_0))
        =\vv.
    \]
    Thus every outgoing local leaf has the same supporting cell $\cc(\suppu{\xevo}(\vv_0))$ and the same oriented straight germ.
    This completes the proof for the outgoing side.
    The incoming side is similar.
\end{proof}

\begin{theorem}\label{thm:X-foli}
    Let $\xevo$ be an arbitrary point in $\CCx(\seeds)$. Then the $\xevo$-\EVO flow induces a piecewise linear $\xevo$-foliation on $\CCx(\seeds)$
    \begin{gather}\label{eq:Xfoli}
        \CCx(\seeds)\backslash\{\xevo,\Xevo\}=\bigcup_{L~\xevo\text{-leaves}} L
    \end{gather}
\end{theorem}

\begin{proof}
By \Cref{lem:sing}, the regular locus is precisely the complement of the two displayed points. 
At any regular point $\ccpt{V}$, \Cref{lem:c-crossing} supplies a unique incoming and a unique outgoing germ. 
Their union with $\ccpt{V}$ is an oriented interval germ: 
\begin{itemize}
    \item in the tangent case $\suppu{\xevo}(\vv)=\suppd{\xevo}(\vv)$ is an ordinary affine integral curve,
    \item while in the crossing case $\suppu{\xevo}(\vv)\ne\suppd{\xevo}(\vv)$ it consists of the two canonical straight germs in their respective cells.
\end{itemize}

Starting with such a germ, continue it through every regular endpoint using \Cref{lem:c-crossing}. 
Equivalently, take the equivalence relation generated by regular local $\xevo$-leaf intervals and their canonical joins at regular boundary points.
The equivalence classes are the maximal leaves. 

Every regular point belongs to a class. 
If two continued curves meet at a regular point, uniqueness of both oriented germs makes them agree there and at every subsequent continuation in either direction. 
They therefore lie in the same maximal class. 
Thus the classes form the disjoint union \eqref{eq:Xfoli}. 

Give each leaf the intrinsic continuation topology obtained by gluing its interval germs. 
There is no branching, and every regular point has an interval neighborhood in this topology. 
Since each segment within a cell is straight, the leaves are piecewise linear in this intrinsic sense.
\end{proof}

In the proof of \Cref{thm:X-foli}, we see that any $\ccpt{V}$ such that $\suppu{\xevo}(\vv)\ne\suppd{\xevo}(\vv)$ connects two non-degenerate local $\xevo$-leaves. Then we can define the following.

\begin{definition}\label{def:c-crossing}
    There is a \emph{cell-crossing} at $\ccpt{V}$ if $\suppu{\xevo}(\vv)\ne\suppd{\xevo}(\vv)$,
    that is, the flow $\flowU{\xevo}$ travels from $\cc(\suppd{\xevo}(\vv))$ to $\cc({\suppu{\xevo}(\vv)})$ / the flow $\flowD{\xevo}$ travels from $\cc(\suppu{\xevo}(\vv))$ to $\cc({\suppd{\xevo}(\vv)})$.
    We call $\vv=\supp(\ccpt{V})$ the \emph{cell-wall} of a cell-crossing at $\ccpt{V}$.
\end{definition}

The cell-crossing on a cell-wall $\vv$ is determined by the parallel flow $\flowdu{\xevo}(\vv)$.
Moreover, it is a continuous refinement of Keller's green mutation as follows.

\begin{lemma}\label{lem:gen-green}
Let ${\xx}$ be any cluster containing $\xevo$.
Then we have the following.
\begin{itemize}
    \item If $u\in\suppu{\xevo}(\vv)\backslash\vv$, then $u$ is green in $\Bgzp{{\xx}}{\vv}$ with respect to ${\xx}$.
    \item If $w\in\suppd{\xevo}(\vv)\backslash\vv$, then $w$ is red in $\Bgzm{{\xx}}{\vv}$ with respect to ${\xx}$.
\end{itemize}
In particular, when $\xevo$ is generic (cf. \eqref{eq:cc}), the upward/downward cell-crossing coincides with the red/green mutation for any cell-wall of co-rank 1.
Therefore, a generic $\xevo$-foliation on $\CCx(\seeds)$ induces $\CEG^\xx(\seeds)$.
\end{lemma}

\begin{proof}
We only show the upward case and the downward case is similar.

If $u\in\suppu{\xevo}(\ccpt{V})\backslash\vv$, then $u\in\supp(\flowu{\xevo}(\vv))$.
Note that $[\flowu{\xevo}(\vv)]_-\in\RR\vv$ by \Cref{def:bong} and $u\not\in\vv$, we have moreover $u\in\supp([\flowu{\xevo}(\vv)]_+)$ and hence
\[
\ubf{e}_x\cvec^{{\xx}}_{u,\Bgzp{{\xx}}{\vv}}>0
\]
for some $x\in\xx$. It implies that
\[
\cvec^{{\xx}}_{u,\Bgzp{{\xx}}{\vv}}\in\ZZp^{n}.
\]
Hence $u$ is green in $\Bgzp{{\xx}}{\vv}$ with respect to ${\xx}$.
\end{proof}

Moreover, the ranks of $\xevo$-supports are invariant under cell-crossing, cf. \Cref{thm:X=xx}.
Hence, we can define the rank of an $\xevo$-leaf as follows.

\begin{definition}\label{def:rankL}
Let $L$ be an $\xevo$-leaf of $\CCx(\seeds)$. Then the \emph{rank} of $L$ is
\begin{gather}
    \rank L\colon=|\suppdu{\xevo}(\supp\ccpt{V})|,\quad\forall~\ccpt{V}\in L\backslash\{\xevo,\Xevo\},
\end{gather}
which does not depend on the choice of $\ccpt{V}$. We say $L$ is \emph{generic}, if $\rank L=\dim\CCx(\seeds)+1$.
\end{definition}

\section{\texorpdfstring{Acyclicity via $\xx$-order}{Acyclicity via x-order}}\label{sec:acyclic}

In this section,
we prove the acyclicity conjecture \Cref{Conj0} for $\xevo$-foliations using the $\xx$-order, defined in terms of the dominant sets introduced in \cite{Cao}.

\subsection{${\xx}$-Order}\label{sec:xorder}\
The acyclic $\xevo$-foliation is defined as follows.

\begin{definition}
Let $\Bar{L}$ be the endpoint completion of $L$.
We say an $\xevo$-leaf $L$ is
\begin{itemize}
    \item \emph{compact}, if $\Bar{L}=L\cup\{\xevo,\Xevo\}\simeq[0,1]$,
    \item \emph{semi-compact}, if $\Bar{L}=L\cup\{\xevo\}\simeq(0,1]$ or $\Bar{L}=L\cup{\Xevo}[0,1)$,
    \item \emph{cyclic}, if $L\simeq\sphere^1$.
\end{itemize}
Moreover, we say an $\xevo$-foliation is
\begin{itemize}
    \item \emph{compact}, if each $\xevo$-leaf is compact,
    \item \emph{semi-compact}, if each $\xevo$-leaf is compact or semi-compact,
    \item \emph{acyclic}, if there are no cyclic $\xevo$-leaves.
\end{itemize}
\end{definition}


For each nonzero polynomial $F(\yy)=\sum_{\bfit{h}}f_{\bfit{h}}\yy^{\bfit{h}}\in\ZZ[\yy]$, let
\begin{gather}\label{eq:tropoly}
    F[-]:\ZZ^n\to\ZZ
\end{gather}
denote the \emph{tropical polynomial} of $F$, defined by
\[
F[\gvec]=\max\{\bfit{h}^T\gvec\mid f_{\bfit{h}}\ne0\}.
\]
We apply this construction to the $F$-polynomials with respect to a cluster ${\xx}$, cf. \Cref{def:Fpoly}.
The following invariant will be central to the argument.

\begin{theorem}\cite[Thm.~4.10]{Cao}\label{prop:tropoly}
For any two cluster variables $v,w$, the number
\[
    \Fpoly^{{\xx}}_w[D\gvec^{{\xx}}_{v}]+\Fpoly^{{\xx}}_v[D\gvec^{{\xx}}_{w}]
\]
does not depend on the choice of the initial cluster ${\xx}$.
\end{theorem}

We next recall the definition of the dominant sets for partial clusters.

\begin{definition}\label{def:dom}
\cite[Def.~7.3, Prop.~7.4(iii)]{Cao}.
The \emph{dominant set} $\dom{{\xx}}$ with respect to ${\xx}$ is an operation
\begin{gather}\label{dom}
    \dom{{\xx}}:\CCx_*(\seeds)\to
    \{\hua{S}\mid\hua{S}\subset\CCx_1(\seeds)\}
\end{gather}
sending a partial cluster $\vv$ to
\begin{align}
    \dom{{\xx}}(\vv)
    &=\{
    w\in\CCx_1(\seeds) \mid \Fpoly^{{\xx}}_w[D\gvec^{{\xx}}_v]=0, \forall~v\in\vv
    \}\\
    &=
    \{
    w\in\CCx_1(\seeds) \mid \Fpoly^{{\xx}}_w[\sum_{v\in\vv} h_v D\gvec^{{\xx}}_v])=0, \forall~h_v\in\ZZ_{>0}
    \}.
\end{align}
\end{definition}

By \cite[Prop.~7.4]{Cao} and the intersection formula in \Cref{def:dom}, the dominant sets satisfy the following property for partial clusters:
\begin{gather}\label{eq:dom}
    {\xx},\vv\subset\dom{{\xx}}(\vv)=\bigcap_{v\in\vv}\dom{{\xx}}(v)
\end{gather}
 Inclusion
of dominant sets, together with reverse inclusion of partial clusters when the dominant sets agree, gives the following partial order on  $\CCx_*(\seeds)$.

\begin{definition}
The \emph{${\xx}$-order} on $\CCx_*(\seeds)$ is defined by
\begin{gather}
    \ww\le_{{\xx}}\uu\iff\dom{\xx}(\uu)\subsetneqq\dom{\xx}(\ww)~\text{or}~(\dom{\xx}(\uu)=\dom{\xx}(\ww)~\&~\uu\subset\ww).
\end{gather}
\end{definition}

\begin{proposition}
The relation $\le_{\xx}$ is a partial order on $\CCx_*(\seeds)$.
\end{proposition}

\begin{proof}
Reflexivity follows from $\vv\subseteq\vv$.  
If $\ww\le_{\xx}\uu$ and $\uu\le_{\xx}\ww$, the dominant sets must agree, since strict inclusion cannot hold in both directions. 
The second clause in the definition then gives 
$\uu\subseteq\ww$ and $\ww\subseteq\uu$; hence $\uu=\ww$.  
For transitivity,
\[
\dom{\xx}(\zz)\subseteq\dom{\xx}(\uu)\subseteq\dom{\xx}(\ww)
\]
whenever $\ww\le_{\xx}\uu\le_{\xx}\zz$.  
If the inclusion between the first and last sets is strict, then $\ww\le_{\xx}\zz$ by the first clause.  
If it is an equality, then all three dominant sets are equal and $\zz\subseteq\uu\subseteq\ww$,
so the second clause applies.
\end{proof}


The restriction of the ${\xx}$-order to $\CCx_n(\seeds)$ agrees with the ${\xx}$-orientation of \Cref{def:green}, in the following sense.

\begin{theorem}\cite[Cor.~7.10, Thm.~7.13]{Cao}\label{thm:xx-order}
For any two clusters $\ww,\uu$, we have
\begin{gather}
    \ww\le_{{\xx}}\uu\impliedby\exists~\text{a sequence of green mutations $\Vec{\mu}$ s.t.}~\ww=\Vec{\mu}(\uu).
\end{gather}
\end{theorem}

Together with \Cref{lem:gen-green}, this implies that all generic $\xevo$-leaves for any $\xevo$-foliation are acyclic .
We now extend this conclusion to arbitrary $\xevo$-leaves.

\subsection{Acyclicity of $\xevo$-foliations}\

To prove the acyclicity of $\xevo$-foliations,
we use the ${\xx}$-order associated to a cluster ${\xx}$ containing $\supp(\xevo)$.

\begin{theorem}\label{thm:X-order}
For any cluster ${\xx}$ containing $\supp(\xevo)$, we have
\begin{gather}\label{eq:u>w>d}
    \dom{{\xx}}(\vv)=\dom{{\xx}}(\suppd{\xevo}(\vv))
    \quad\text{and}\quad
    \dom{{\xx}}(\suppu{\xevo}(\vv))\subsetneqq\dom{{\xx}}(\vv)
\end{gather}
for any $\suppu{\xevo}(\vv)\ne\suppd{\xevo}(\vv)$.
\end{theorem}

\begin{proof}
For the equality in \eqref{eq:u>w>d}, we adapt the proof of \cite[Thm.~7.8]{Cao}.
Recall the natural embedding \eqref{eq:gf-embedding} associated to ${\xx}$ containing $\supp(\xevo)$.
Consider the integer vector
\[
\bfit{g}_0=\sum_{v\in\vv} \alpha_v \gvec^{{\xx}}_{v}\in\sigma^{{\xx}}(\vv)
\]
with $\alpha_v\in\ZZ_{>0}$.
For each $w\in\Bgzm{\xx}{\vv}$, put
\[
    \bfit{n}_w
    :=
    d_w^{-1}D
    \cvec^{\xx}_{
        w,\Bgzm{\xx}{\vv}}.
\]
By the tropical duality \eqref{eq:dual2}, these vectors satisfy
\[
    \bfit{n}_w^{T}\gvec_z^{\xx}
    =\delta_{wz},
    \qquad
    w,z\in\Bgzm{\xx}{\vv}.
\]
In particular, they are the inward normal vectors dual to the generators of the $g$-cone $\con^{\xx}(\Bgzm{\xx}{\vv})$.
By the sign condition of Bongartz co-completion,
\[
    \bfit{n}_w\leq 0,
    \qquad
    \forall w\in
    \Bgzm{\xx}{\vv}
    \setminus\vv.
\]
Set
\[
    \bfit{p}
    :=
    \sum_{
        w\in
        \Bgzm{\xx}{\vv}
        \setminus\vv}
    \bfit{n}_w.
\]
For every
$w\in
\Bgzm{\xx}{\vv}
\setminus\vv$,
the coordinatewise nonpositivity of the summands gives
\[
    \bfit{n}_w^{T}\bfit{p}
    =
    \sum_{
        z\in
        \Bgzm{\xx}{\vv}
        \setminus\vv}
    \bfit{n}_w^{T}\bfit{n}_z
    \geq
    \bfit{n}_w^{T}\bfit{n}_w
    >0.
\]
Moreover,
\[
    \bfit{n}_w^{T}\gvec_0
    =
    \begin{cases}
        \alpha_w, & w\in\vv,\\
        0, & w\in
        \Bgzm{\xx}{\vv}
        \setminus\vv.
    \end{cases}
\]
Since $\alpha_w>0$ for $w\in\vv$, we may choose a sufficiently small rational number $\varepsilon>0$ such that
\[
    \gvec:=\gvec_0+\varepsilon\bfit{p}
    \quad\text{satisfies}\quad
    \bfit{n}_w^{T}\gvec>0
    \quad
    \text{for every }
    w\in\Bgzm{\xx}{\vv}.
\]
Consequently, 
\[
    \gvec\in\sigma^{\xx}\bigl(\Bgzm{\xx}{\vv}\bigr)^\circ.
\]

Taking $\alpha_w=\varepsilon$ for $w\in\Bgzm{\xx}{\vv}\setminus\vv$, we obtain
\[
    \gvec
    =
    \gvec_0+
    \sum_{
        w\in
        \Bgzm{\xx}{\vv}
        \setminus\vv}
    \alpha_w d_w^{-1}
    \bigl(
        D\cvec^{\xx}_{
            w,\Bgzm{\xx}{\vv}}
    \bigr).
\]
All vectors and coefficients involved are rational.
Clearing denominators, choose $N\in\mathbb Z_{>0}$ such that
$N\alpha_wd_w^{-1}\in\mathbb Z_{>0}$ for all
$w\notin\vv$ and
\[
    N\gvec
    =
    \sum_{
        v\in\Bgzm{\xx}{\vv}}
    \alpha_v^*\gvec_v^{\xx},
    \qquad
    \alpha_v^*\in\mathbb Z_{>0}.
\]
Thus \Cref{def:bong} gives
\begin{gather}\label{eq:diff>0}
    N\bfit{g}_0-N\bfit{g}=
-\sum_{w\in\Bgzm{{\xx}}{\vv}\backslash\vv} N \alpha_w d^{-1}_w(D\cvec^{{\xx}}_{w,\Bgzm{\xx}{\vv}})\in\ZZp^n.
\end{gather}
For any $z\in\dom{{\xx}}(\vv)$, we have $\Fpoly^{{\xx}}_z[DN\gvec_0]=0$; thus $\bfit{h}^T DN\gvec_0\le0$ for every exponent vector $\bfit{h}\in\ZZp^{n}$ occurring in $\Fpoly^{{\xx}}_z$.
Using \eqref{eq:diff>0}, we obtain
\[
    0\le
    \Fpoly^{{\xx}}_z[DN\bfit{g}]=\max(\bfit{h}^T\cdot DN\gvec)=\max(\bfit{h}^T\cdot DN\gvec_0-\bfit{h}^T\cdot(DN\bfit{g}_0-DN\bfit{g}))\le 0.
\]
Here the first inequality follows from the constant term $1$ of every
$F$-polynomial.
Hence $F_z^{\xx}[DN\gvec]=0$.
Since every coefficient $\alpha_v^*$ is positive,
the support-independence of dominant sets
\cite[Proposition~7.4(iii)]{Cao} gives
\[
    z\in
    \dom{\xx}
    \bigl(\Bgzm{\xx}{\vv}\bigr).
\]
Thus
\[
    \dom{\xx}(\vv)
    \subseteq
    \dom{\xx}
    \bigl(\Bgzm{\xx}{\vv}\bigr).
\]
On the other hand,
\[
    \vv
    \subseteq
    \suppd{\xevo}(\vv)
    \subseteq
    \Bgzm{\xx}{\vv},
\]
so \Cref{def:dom} yields
\[
    \dom{\xx}
    \bigl(\Bgzm{\xx}{\vv}\bigr)
    \subseteq
    \dom{\xx}
    \bigl(\suppd{\xevo}(\vv)\bigr)
    \subseteq
    \dom{\xx}(\vv).
\]
Therefore all three dominant sets coincide, proving the equality in \eqref{eq:u>w>d}.

For the strict inclusion, put $\uu=\Bgzp{\xx}{\vv}$ and choose
\[
u\in\suppu{\xevo}(\vv)\setminus\vv,
\qquad w=\DT{(u)}_\vv\in\Bgzm{\xx}{\vv},
\]
using the bijection \eqref{eq:Bgz-DT}.  
In the reduced cluster pattern at $\vv$, with initial cluster $\uu\setminus\vv$, the corresponding DT variable satisfies
$\gvec^{\uu}_w=-\bfit e_u$, whereas $\gvec^{\uu}_u=\bfit e_u$.
Thus $u$ and $w$ are incompatible: if they belonged to a common cluster, their
$g$-vectors would be part of a basis and hence linearly independent.

By the compatibility criterion for the $F$-invariant
\cite[Thm.~4.19]{Cao}, their $F$-invariant is positive.  Since
$\Fpoly^{\uu}_u=1$, this says
\begin{gather}\label{eq:F3}
\Fpoly^{\uu}_w[D\gvec^{\uu}_u]
+\Fpoly^{\uu}_u[D\gvec^{\uu}_w]
=\Fpoly^{\uu}_w[D\bfit e_u]>0.
\end{gather}
On the other hand, $u\in\uu\subset\dom{\xx}(\vv)$, and the equality established above gives 
$\dom{\xx}(\vv)=\dom{\xx}(\Bgzm{\xx}{\vv})$.  Since
$w\in\Bgzm{\xx}{\vv}$, it follows that
\[
\Fpoly^{\xx}_u[D\gvec^{\xx}_w]=0.
\]
Applying the invariance in \Cref{prop:tropoly} to $u,w$ now yields
\[
\Fpoly^{\xx}_w[D\gvec^{\xx}_u]
=\Fpoly^{\uu}_w[D\gvec^{\uu}_u]>0.
\]
Therefore $w\notin\dom{\xx}(\suppu{\xevo}(\vv))$, whereas
$w\in\dom{\xx}(\vv)$.  The inclusion is strict.
\end{proof}

We now apply \Cref{thm:X-order} to prove the acyclicity conjecture for $\xevo$-foliations.

\begin{theorem}\label{thm:acyclic}
    Any $\xevo$-foliation is acyclic.
\end{theorem}

\begin{proof}
Suppose there exists a cyclic $\xevo$-leaf $L$.
By \Cref{thm:X-order}, we have
\begin{gather}\label{eq:gXorder}
    \suppd{\xevo}(\vv)\le_{{\xx}}\vv\le_{{\xx}}\suppu{\xevo}(\vv)
\end{gather}
for cell-crossing at $\vv=\supp(\ccpt{V})$ and any cluster ${\xx}$ containing $\supp(\xevo)$.
Because a cyclic leaf is compact, there are  finitely many such cell-crossing.
Following the oriented leaf once around therefore gives a finite cycle
\[
\vv_1\le_{{\xx}}\vv_2\le_{{\xx}}\ldots\le_{{\xx}}\vv_k\le_{{\xx}}\vv_1
\]
in $(\CCx_*(\seeds),\le_{{\xx}})$.
Since $\le_{{\xx}}$ is a partial order, we have $\vv_1=\ldots=\vv_k$ and
$L\subset\cc(\vv_1)$.
This contradicts the linearity of local $\xevo$-foliation \eqref{eq:local-foli}.
\end{proof}


\begin{remark}
    The partial order generated by \eqref{eq:gXorder} need not recover the ${\xx}$-order even if $\xevo\in\cc(\xx)^\circ$ is generic.
\end{remark}


For cluster complexes of finite type, \Cref{thm:acyclic} implies that every $\xevo$-foliation is compact. Moreover, a generic
$\xevo$-foliation induces a homotopy equivalence between the cluster complex and a sphere.

\begin{corollary}
If $|\CCx_*(\seeds)|<\infty$, then any $\xevo$-foliation on $\CCx(\seeds)$ is compact.
\end{corollary}

\begin{proof}
If $|\CCx_*(\seeds)|<\infty$, then each $\xevo$-leaf $L$ admits finitely many cell-crossings.
Indeed, at every cell-crossing the order changes strictly by \Cref{thm:X-order}, so a leaf cannot pass through the same finite collection of cells indefinitely.
Since the finite cluster complex is compact, every end of a non-cyclic leaf has an accumulation point.  The local extension property of
the global foliation shows that such an end continues unless the accumulation
point is one of the two singularities $\xevo,\Xevo$.  Thus every leaf is either
compact or cyclic.  By \Cref{thm:acyclic}, every leaf is compact, so the
$\xevo$-foliation is compact.
\end{proof}

\begin{corollary}
If $|\CCx_*(\seeds)|<\infty$, then $\CCx(\seeds)\simeq\sphere^{\dim\CCx(\seeds)}$.
\end{corollary}

\begin{proof}
The hypothesis says that the cluster pattern is of finite type.  
Choose $\xevo$ to be a generic point, the resulting $\xevo$-foliation induces a homotopy equivalence
\begin{gather}
    \CCx(\seeds)\simeq\susp\partial(\cc(\supp(\xevo)))\simeq\susp \sphere^{\dim\CCx(\seeds)-1}\simeq\sphere^{\dim\CCx(\seeds)}.
\end{gather}
\end{proof}



\Cref{thm:X-order} also provides a criterion for one-sided compactness of $\xevo$-leaves.

\begin{corollary}
Let $L$ be an $\xevo$-leaf, let $\xx$ be
a cluster containing $\operatorname{supp}(\xevo)$, and let
$\vv$ be a partial cluster such that
\[
   \cc(\vv)^\circ
    \cap L\neq\varnothing.
\]
Then the following statements hold.
\begin{enumerate}
    \item If $|\dom{\xx}(\vv)|<\infty$, 
    then the upward end of $L$ is $\xevo$.
    \item If
    $|\CCx_1(\seeds)\setminus\dom{\xx}(\vv)|<\infty$,
    then the downward end of $L$ is
    $\Xevo$.
\end{enumerate}
In either case, $L$ is compact or semi-compact.
If both finiteness conditions hold, then $L$ is compact.
\end{corollary}

\begin{proof}
Choose $\ccpt{V}\in\cc(\vv)^\circ\cap L$.
By \Cref{thm:X-order}, at each cell-crossing with cell-wall $\ww$, we have
\[
    \dom{\xx}
        \bigl(\suppu{\mathrm{X}}(\ww)\bigr)
    \subsetneq
    \dom{\xx}(\ww)
    =
    \dom{\xx}
        \bigl(\suppd{\mathrm{X}}(\ww)\bigr).
\]
Thus, the dominant sets strictly decrease at each cell-crossing along the upward flow. 
If $\dom{\xx}(\vv)$ is finite, there can be only finitely many such crossings after $\ccpt{V}$.

Similarly, the dominant sets strictly increase along the downward flow, so their complements in $\CCx_1(\seeds)$ strictly decrease. 
Hence the second finiteness condition implies that there are only finitely many
downward cell-crossings after $\ccpt{V}$.

The orientation determines which singularity occurs.
The upward end is $\xevo$ and the downward end is $\Xevo$.

\end{proof}

\appendix

\section{Categorical settings}\label{app:cats}
In this appendix, we compare our setting (cluster algebras) with two categorical settings:
the 2-Calabi-Yau categories studied in the prequel and the $\tau$-tilting theory.
\subsection{X-\EVO flows for 2-Calabi--Yau cluster categories}\label{sec:2CY}\

\paragraph{\textbf{cluster tilting set and cluster complex}}

We recall the categorical construction of \cite{QT} and explain its relation with the construction in this paper.
Let $\C$ be a Hom-finite, Krull-Schmidt, $\k$-linear triangulated category,
where $\k$ is algebraically closed. 
Assume that $\C$ is $2$-Calabi-Yau,
i.e.
\[
    \Hom_{\C}(A,B)\cong
    D\Hom_{\C}(B,A[2])
\]
functorially in $A,B$, and that $\C$ admits cluster tilting objects.

\begin{definition}[{\cite[Def.~2.2,~Rem.~2.3~and~Def.~2.5]{QT}}]\label{def:app-cluster}
A basic object $\pcto{V}$ is \emph{cluster tilting} if
\[
    \Add\pcto{V}
    =\{A\in\C:\Ext^1_{\C}(\pcto{V},A)=0\}.
\]
Its indecomposable summands form a cluster tilting set $\pcts{V}$.
A \emph{partial cluster tilting set} is a subset of a cluster tilting set.
The \emph{cluster complex} $\CCx(\C)$ has these partial sets as simplices;
compatibility is given by the vanishing of $\Ext^1_{\C}$.
\end{definition}

In particular, we need the following mild condition on $\C$ for simplicity.
\begin{itemize}
  \item[$\bigstar$] Any rigid objects in $\C$ is a partial cluster.
\end{itemize}

For comparison, assume $\C$ gives a cluster categorification of a skew-symmetric cluster pattern $\seeds(\C)$ satisfying \Cref{ass0}. 
We work with the reachable cluster tilting sets. 
Accordingly, let $\CCx^{\circ}(\C,\pcts{X})$ denotes the mutation-reachable complex with respect to a cluster $\pcts{X}$, and then
\begin{equation}\label{eq:app-complex}
    \CCx^\circ(\C,\pcts{X})\cong\CCx(\seeds(\C)).
\end{equation}
This qualification is unnecessary when all cluster tilting sets under consideration are mutation-reachable.

\paragraph{\textbf{Evolving triangles and partial Bongartz (co-)completions.}}\

Fix an indecomposable rigid object $X$ and a partial cluster tilting set $\pcts{V}$. 
We recall the definition of evolving triangles.

\begin{definition}[{\cite[Def./Lem.~3.2, Def.~3.5]{QT}}]
\label{def:app-triangle}
The \emph{downward $X$-evolving triangle} for $\pcts{V}$ is
\begin{equation}\label{eq:app-triangle}
    X\xrightarrow{g}W\longrightarrow U\xrightarrow{f}X[1],
\end{equation}
where $f$ is a minimal right $\Add\pcto{V}$-approximation of $X[1]$.
The \emph{upward $X$-evolving triangle} is defined by requiring $g$
to be a minimal left $\Add\pcto{V}$-approximation of $X$ instead.
\end{definition}

These triangles exist and are unique up to isomorphism. 
If $\pcts{V}$ is a cluster tilting set, the two triangles coincide.
In either case, $\Add W\cap\Add U=0$ and
$W\oplus \pcto{V}\oplus U$ is a partial cluster in which the $\xevo$-evo flow supports \cite[Lem.~3.7]{QT}. 

Hence the evolving triangles provide the categorical replacement for the partial Bongartz (co-)completions.
Let $\pcts{V}$ be a partial cluster.
Note that $U\in\Add\pcto{V}$ in the downward
case and $W\in\Add\pcto{V}$ in the upward case \cite[Rem.~3.6]{QT}.
\begin{itemize}
    \item If \eqref{eq:app-triangle} is downward for $\pcts{V}$, we have
    \begin{equation}\label{eq:app-bongartz-down}
        \Bgzm{\{X[1]\}}{\pcts{V}}
        =\Ind\Add(W\oplus\pcto{V}\oplus U)
        =\Ind\Add(W\oplus\pcto{V}).
    \end{equation}
    \item If \eqref{eq:app-triangle} is upward for $\pcts{V}$, we have
    \begin{equation}\label{eq:app-bongartz-up}
        \Bgzp{\{X[1]\}}{\pcts{V}}
        =\Ind\Add(W\oplus\pcto{V}\oplus U)
        =\Ind\Add(\pcto{V}\oplus U).
    \end{equation}
\end{itemize}
The partial Bongartz (co-)completion with respect to a partial cluster $\pcts{X}[1]$ is given by union, cf. \Cref{rem:union}.
\begin{gather}
    \Bgzpm{\pcts{X}[1]}{\pcts{V}}=\bigcup_{X\in\pcts{X}}\Bgzpm{\{X[1]\}}{\pcts{V}}.
\end{gather}

\paragraph{\textbf{$\xevo$-evo flow and $\xevo$-foliation.}}\

For a nonzero rigid object $A=\bigoplus_j A_j^{\oplus m_j}$, with
pairwise non-isomorphic indecomposable summands, put
\[
    \real{A}:=\frac{\sum_j m_j A_j}{\sum_j m_j}\in\CCx(\C),
    \qquad \real{0}:=0.
\]
Here $\real{0}$ is a zero vector, rather than a point of the complex.

\begin{definition}[{\cite[Def.~3.9 and Def.~3.11]{QT}}]\label{def:app-flow}
For $\ccpt{V}\in\cc(\pcts{V})^\circ$, use the downward triangle
\eqref{eq:app-triangle} to define
\begin{equation}\label{eq:app-down}
    \Dflow{X}(\ccpt{V})=
    \begin{cases}
        X-\ccpt{V},&U=0,\\
        \real{W}-\real{U},&W\ne0\ne U,\\
        \ccpt{V}-X[1],&W=0.
    \end{cases}
\end{equation}
Using the upward triangle instead, define
\begin{equation}\label{eq:app-up}
    \Uflow{X}(\ccpt{V})=
    \begin{cases}
        \ccpt{V}-X,&U=0,\\
        \real{U}-\real{W},&W\ne0\ne U,\\
        X[1]-\ccpt{V},&W=0.
    \end{cases}
\end{equation}
For $\xevo=\sum_i\coex_iX_i\in\CCx(\C)$, extend these definitions by
\[
    \DUflow{\xevo}=\sum_i\coex_i\DUflow{X_i}
    \quad\text{and}\quad \xevo[1]=\sum_i\coex_iX_i[1].
\]
\end{definition}

Take $\pcts{V}$ to be the Bongartz completion of $\supp\ccpt{V}$ with respect to a cluster containing $X[1]$.
Use the exchange-matrix convention in which the transposed $g$-vector
$\tgvec_A^{\pcts{V}}$ is the coordinate vector of the index $\on{ind}_{\pcts{V}}A$.
Rotation of \eqref{eq:app-triangle} gives
\begin{equation}\label{eq:app-index}
    \on{ind}_{\pcts{V}}X[1]=[U]-[W]
    \quad\text{and}\quad
    \tgvec_{X[1]}^{\pcts{V}}=(u_j-w_j)_j,
\end{equation}
where $U=\bigoplus_jV_j^{\oplus u_j}$ and
$W=\bigoplus_jV_j^{\oplus w_j}$. 
Minimality implies $u_jw_j=0$.
Consequently, 
\begin{gather}
    \normal_+(\tgvec^\pcts{V}_{\source})=\real{U}
    \quad\text{and}\quad
    \normal_-(\tgvec^\pcts{V}_{\source})=\real{W}.
\end{gather}
Thus the nontrivial triangles give the parallel part, while the two trivial triangles give the radial terms. 
Subtracting the coefficient sum times $\ccpt{V}$ gives exactly the radial modification in \eqref{eq:app-down} and \eqref{eq:app-up}.
Extending linearly in $\xevo$,
\Cref{def:app-flow} coincide with \Cref{def:Evo} in the sense that
\begin{equation}\label{eq:app-fields}
    \Dflow{\xevo}=\flowD{\xevo[1]}
    \quad\text{and}\quad
    \Uflow{\xevo}=\flowU{\xevo[1]}.
\end{equation}
In particular, the DT convention of this paper gives
\begin{equation}\label{eq:app-shift}
    \DT{X}=X[-1],\qquad \iDT{X}=X[1],
    \qquad
    \flowDU{\xevo}=\DUflow{\Xevo}.
\end{equation}

In the categorical setting there, we also proved that $\xevo$-evo flows induced $\xevo$-foliations (\cite[Thm.~5.10]{QT}), 
as we do here in \Cref{thm:X-foli}.

\begin{remark}
    In \cite{HZ}, they prove the acyclicity conjecture \cite[Conj.~1.7]{QT} in the category setting, as we do in \Cref{thm:acyclic} in the cluster-algebraic setting.  
    Moreover, they also prove \cite[Conj.~1.8]{QT}, that is, semi-compactness of any $\xevo$-foliation in the cluster complex of Euclidean type.
\end{remark}

\paragraph{\textbf{The correspondence in a nutshell.}}\

The following dictionary uses the categorical parameter $\xevo$;
the cluster-algebraic parameter in the main body is $\xevo[1]$.
\begin{center}
\small
\renewcommand{\arraystretch}{1.18}
\begin{tabular}{@{}p{.47\linewidth}p{.47\linewidth}@{}}
\hline
Categorical setting & Cluster-algebraic setting\\
\hline
Rigid indecomposables
& Cluster variables\\
(Partial) cluster tilting sets
& (Partial) clusters\\
Minimal right/left approximation triangles
& Partial Bongartz co-completion/completion\\
$[U]-[W]=\on{ind}_{\pcts{V}}X[1]$
& Transposed $g$-vector $\tgvec_{X[1]}^{\vv}$\\
The parallel part $\real{W}-\real{U}$
& Downward parallel part $\flowd{X[1]}$\\
The $\xevo$-evo flow $\DUflow{\xevo}$ 
&The $\xevo$-evo flow $\flowDU{\xevo[1]}$\\
Shift $[-1]$
& $\DT{(-)}$\\
\hline
\end{tabular}
\end{center}

\subsection{Partial Bongartz (co)-completions in $\tau$-tilting theory}
\label{app:tau-tilting}\

\paragraph{\textbf{The support $\tau$-tilting pairs}}\
We recall the module-theoretic interpretation of the constructions in
\Cref{sec:suppBgz}. 
Throughout, we denote by
\begin{itemize}
    \item $\pcts{X}$ a basic finite-dimensional $\k$-algebra with $n$ simple modules,
    \item $\on{mod}\pcts{X}$ the category of finite-dimensional right $\pcts{X}$-modules,
    \item $\tau$ the Auslander-Reiten translation,
    \item $|\pcto{V}|$ the number of non-isomorphic indecomposable summands of a finite dimensional right $\pcts{X}$-module $\pcto{V}$ and 
    \item $\on{Fac}\pcto{V}$ the quotients of objects in $\Add \pcto{V}$.
    \item $(-)^\perp$ and ${^\perp}(-)$ the  perpendiculars with respect to $\Hom_\pcts{X}$.
\end{itemize}
The key ingredient in the $\tau$-tilting theory is the \emph{support $\tau$-tilting pairs}, since they bridges various subjects as follows.
\begin{itemize}
    \item $\on{s\tau-tilt}\pcts{X}$ is the set of the support $\tau$-tilting pairs,
    \item $\on{f-tors}\pcts{X}$ is the set of \emph{functorially finite torsion classes} and
    \item $\on{2-silt}\pcts{X}$ is the set of \emph{2-term silting complexes} in $\per\pcts{X}$.
\end{itemize}

We first recall the definition of the support $\tau$-tilting pairs.

\begin{definition}[{\cite[\S2]{AIR}}]\label{def:app-tau-pair}
    A pair $(\pcto{V},P)$ in $\on{mod}\pcts{X}$ with $P$ projective, is \emph{$\tau$-rigid} if
    \[
     \Hom_\pcts{X}(\pcto{V},\tau \pcto{V})=0,\qquad \Hom_\pcts{X}(P,\pcto{V})=0.
    \]
    Moreover, it is \emph{support $\tau$-tilting} if $|\pcto{V}|+|P|=n$ and \emph{almost complete} if $|\pcto{V}|+|P|=n-1$.
    A support $\tau$-tilting pair $(\pcto{V},P)$ is determined by $\pcto{V}$ uniquely, cf. \cite[Prop.~2.3]{AIR}.
\end{definition}


\paragraph{\textbf{2-term silting and g-fans.}}\

The first correspondence is between the support $\tau$-tilting pairs and the 2-term silting objects in $\per \pcts{X}$.

\begin{definition}
A complex $\mathbf{P}_\bullet=(P_i,d_i)$ in $\per\pcts{X}$ is 
\begin{itemize}
    \item \emph{2-term} if $P_i=0$ for $i\ne0,1$,
    \item \emph{presilting} if $\Hom(\mathbf{P}_\bullet,\mathbf{P}_\bullet[i])=0$ for any $i>0$ and
    \item \emph{silting} if it is presilting and $\on{thick}\mathbf{P}_\bullet=\per \pcts{X}$. 
\end{itemize}
\end{definition}

\begin{theorem}[{\cite[Thms.~3.2]{AIR}}]
\label{thm:app-tau-dictionary}
There is a bijection
\[
 \on{s\tau-tilt}\pcts{X}\ \longleftrightarrow\
 \on{2-silt}\pcts{X}
\]
given by
\begin{align}
    \on{s\tau-tilt}\pcts{X}\ni (\pcto{V},P)&\mapsto C_{(\pcto{V},P)}=\bigl(P_1^\pcts{V}\oplus P\xrightarrow{(d_V\;0)}P_0^\pcts{V}\bigr)\in\on{2-silt}\pcts{X}\\
    \on{2-silt}\pcts{X}\ni \mathbf{P}&\mapsto H_0(\mathbf{P})\in\on{s\tau-tilt}\pcts{X},
\end{align}
where $P_1^\pcts{V}\xrightarrow{d_V}P_0^\pcts{V}\to V\to0$ is a minimal projective presentation. 
The bijection can be restricted to $\tau$-rigid pairs, identifying them with two-term presilting complexes. 
\end{theorem}




\paragraph{\textbf{Torsion class and Bongartz completion}}\

The Bongartz endpoint/completion is an important concept in the $\tau$-tilting theory, obtained by the correspondence between the support $\tau$-tilting pairs and the functorially finite torsion classes.
Cao-Gyoda-Yurikusa realize it via the $c$-vectors in cluster algebra, cf. \cite{CGY}.  

\begin{definition}[{\cite[\S2]{AIR}}]
A \emph{torsion class} is a full subcategory of $\on{mod}\pcts{X}$ closed under quotients, extensions and finite direct sums. 
It is \emph{functorially finite} if every module admits both a left and a right approximation by its objects.
\end{definition}

The $\tau$-tilting theory constructs the following correspondences.

\begin{theorem}[{\cite[Thms.~2.7]{AIR}}]
\label{thm:app-tau-dictionary2}
There is a bijection
\[
 \on{f-tors}\pcts{X}\ \longleftrightarrow\
 \on{s\tau-tilt}\pcts{X}.
\]
given by
\begin{gather}
    \on{f-tors}\pcts{X}\ni\hua{T}\mapsto P(\hua{T})\in\on{s\tau-tilt}\pcts{X}
    \quad\text{and}\quad
    \on{s\tau-tilt}\pcts{X}\ni \pcto{V}\mapsto \on{Fac}\pcto{V}\in\on{f-tors}\pcts{X},
\end{gather}
where $P(\mathcal T)$ be the basic direct sum of its indecomposable Ext-projectives, which is support $\tau$-tilting.
\end{theorem}

The definition of the Bongartz (co-)completions arises from a weaker version of the theorem above,
where a $\tau$-rigid pair is associated with an interval in $\on{f-tors}\pcts{X}$ rather than a torsion class.

\begin{theorem}[{\cite[Prop.~2.9, Thm.~2.10]{AIR};
\cite[Thm.~2.8]{AIR24}}]\label{thm:app-tau-endpoints}
Let $(\pcto{V},P)$ be $\tau$-rigid. 
The support $\tau$-tilting pairs containing $(\pcto{V},P)$ correspond exactly to an interval in $\on{f-tors}A$ as follows.
\begin{equation}\label{eq:app-tau-interval}
 [\on{Fac}\pcto{V},{}^\perp(\tau \pcto{V})\cap P^\perp]=
 \bigl\{\mathcal T\in\on{f-tors}\pcts{X}:
 \on{Fac}\pcto{V}\subseteq\mathcal T
 \subseteq{}^\perp(\tau \pcto{V})\cap P^\perp\bigr\}.
\end{equation}
Both endpoints are functorially finite.
We say
\begin{itemize}
    \item the upper endpoint gives $\Bgzp{\pcts{X}}{\pcts{V}}$, the \emph{Bongartz completion} with respect to $\pcts{X}$ and
    \item the lower endpoint gives $\Bgzm{\pcts{X}}{\pcts{V}}$, the \emph{Bongartz co-completion} with respect to $\pcts{X}$.
\end{itemize}
Thus every $\tau$-rigid pair has a completion.
\end{theorem}
\paragraph{\textbf{Partial Bongartz completion in $\on{mod}\pcts{X}$}}\

Inspired by the $\xx$-support and its relation with evolving triangles, 
the partial Bongartz (co-)completion $\Bgzpm{X_i}{\pcts{V}}$ with respect to $X_i\in\pcts{X}$ can be given by approximations in $\on{mod}\pcts{X}$.
Denote by 
\[
\hua{T}^\uparrow_\pcts{V}={}^\perp(\tau \pcto{V})\cap P^\perp
\quad\text{and}\qquad
\hua{T}^\downarrow_\pcts{V}=\on{Fac}\pcto{V}
\]
the torsion classes corresponding to $\Bgzpm{\pcts{X}}{\pcts{V}}$.
For the Bongartz completion, 
taking the minimal left $\hua{T}^\uparrow_{\pcts{V}}$-approximation for $X_i$
\[
    f'_i\colon X_i\longrightarrow W^\uparrow_i
\]
Then we have
\begin{gather}
    \Bgzp{\{X_i\}}{\pcts{V}}=(\on{basic}(\pcto{V}\oplus W^\uparrow_i), P).    
\end{gather}

For the Bongartz co-completion,
taking the minimal left $\Add\pcto{V}$-approximation
\[
 f_i\colon X_i\longrightarrow W^\downarrow_i.
\]
Set $U^\downarrow_i=\on{Coker} f_i$. We have
\[
    X_i\xrightarrow{f_i}W^\downarrow_i\longrightarrow U^\downarrow_i\longrightarrow0,
    \qquad
    W^\downarrow_i\in\Add\pcto{V},\quad
    U^\downarrow_i\in\Add P(\mathcal{T}_{\pcts{V}}^{\downarrow}),
\]
where the last assertion follows from \cite[Lem.~2.20]{AIR}.
To retain the projective support part, choose a minimal projective presentation
\[
    P_1^{W^\downarrow_i}\xrightarrow{d_i}P_0^{W^\downarrow_i}
        \xrightarrow{p_i}W^\downarrow_i\longrightarrow0
\]
and a lift $\widetilde f_i\colon X_i\to P_0^{W^\downarrow_i}$ of $f_i$.
Define the projective module $Q_i$ by the decomposition of projective
presentations
\[
 \bigl(P_1^{W^\downarrow_i}\oplus X_i
       \xrightarrow{(d_i\ \widetilde f_i)}P_0^{W^\downarrow_i}\bigr)
 \simeq
 \bigl(P_1^{U^\downarrow_i}\xrightarrow{d_{U^\downarrow_i}}P_0^{U^\downarrow_i}\bigr)
 \oplus(Q_i\longrightarrow0)
 \oplus(L_i\xrightarrow{1}L_i),
\]
where the first summand is a minimal projective presentation of $U_i$ and $L_i$ are the redundant identity summands.
Hence
\begin{gather}
    \Bgzm{\{X_i\}}{\pcts{V}}=(\on{basic}(\pcto{V}\oplus U^\downarrow_i),\on{basic}(P\oplus Q_i)).
\end{gather}

\paragraph{The correspondence in a nutshell.}
\begin{center}
\small
\renewcommand{\arraystretch}{1.15}
\begin{tabular}{@{}p{.45\linewidth}p{.49\linewidth}@{}}
\hline
Categorical/cluster setting & $\tau$-tilting setting\\
\hline
partial clusters; clusters 
 & $\tau$-rigid pairs; support $\tau$-tilting pairs\\
Summands in $\Add T[1]$
 & Projective support part $(0,P)$\\
transposed $g$-vectors
 & g-Classes of two-term silting complexes\\
Clusters containing a fixed partial cluster
 & The torsion-class interval \eqref{eq:app-tau-interval}\\
$\Bgzp{\ol{\xx}}{\vv}$/$\Bgzm{\ol{\xx}}{\vv}$
 & Upper/lower endpoint\\
Partial Bongartz (co-)completions $\Bgzpm{\xx}{\vv}$
 & Common lifted summands as $\ol{\xx}\supset\xx$ varies\\
\hline
\end{tabular}
\end{center}

\section{\texorpdfstring{$A_3$-Example}{A3-Example}}

Consider the cluster pattern $\seeds_{\DynA{3}}$ of type $\DynA{3}$ generated by the initial seed
\begin{gather}
    {\xx}=(x_1,x_2,x_3),~\yy=(y_1,y_2,y_3)~~\text{and}~~
    B=\left(
    \begin{array}{cccc}
        0 & -1 & -1\\
        1 & 0 & 0\\
        1 & 0 & 0
    \end{array}
    \right)
\end{gather}
The remaining six cluster variables $\{x_4,\ldots,x_9\}$ can be expressed in terms of $\ol{\xx}$, $\yy$, and $\hat{\yy}$ (see \eqref{eq:haty} for the definition of the $\hat{y}$-variables) as follows:
\begin{align*}
    &x_4=x_1^{-1}\cdot\frac{1+\hat{y}_1}{1\oplus y_1},\quad
    x_5=x_3^{-1}\cdot\frac{1+\hat{y}_3+\hat{y}_1\hat{y}_3}{1\oplus y_3\oplus y_1 y_3},\quad
    x_6=x_2^{-1}\cdot\frac{1+\hat{y}_2+\hat{y}_1\hat{y}_2}{1\oplus y_2\oplus y_1 y_2},\\
    &x_7=x_1x_2^{-1}x_3^{-1}\cdot\frac{1+\hat{y}_2+\hat{y}_3+\hat{y}_2\hat{y}_3+\hat{y}_1\hat{y}_2\hat{y}_3}{1\oplus y_2\oplus y_3\oplus y_2y_3\oplus y_1y_2y_3},\\
    &x_8=x_1x_3^{-1}\cdot\frac{1+\hat{y}_3}{1\oplus y_3},\quad
    x_9=x_1x_2^{-1}\cdot\frac{1+\hat{y}_2}{1\oplus y_2}.
\end{align*}
The corresponding $F$-polynomials follow immediately from the first formula in \Cref{thm:CGF}.
For instance, the $F$-polynomial of $x_7$ is
\[
\Fpoly^\xx_{x_7}(\yy)=1+y_2+y_3+y_2y_3+y_1y_2y_3.
\]

For any $\vv$ is a partial cluster but not a cluster,
we compute their Bongartz (co)-completions with respect to ${\xx}$ in \Cref{tb:A3-Bgz1}.

\begin{table}[h]\label{tb:A3-Bgz1}
\begin{tabular}{|c||c|c|}
\hline
$\Bgzpm{{\xx}}{\vv}$ & $\vv\to\Bgzp{{\xx}}{\vv}$ & $\vv\to\Bgzm{{\xx}}{\vv}$ \\
\hline\hline
$(x_1,x_2,x_3)$  & $x_1,x_2,x_3,(x_1,x_2),(x_2,x_3),(x_1,x_3)$ & $\empty$\\
\hline
$(x_2,x_3,x_4)$ & $x_4,(x_2,x_4),(x_3,x_4)$ & $(x_2,x_3)$\\
\hline
$(x_2,x_4,x_5)$ & $(x_4,x_5)$ & $x_2,(x_2,x_4),(x_2,x_5)$\\
\hline
$(x_3,x_4,x_6)$ & $(x_4,x_6)$ & $x_3,(x_3,x_4),(x_3,x_6)$\\
\hline
$(x_4,x_5,x_6)$ & $\empty$ & $x_4,x_5,x_6,(x_4,x_5),(x_5,x_6),(x_4,x_6)$\\
\hline
$(x_5,x_6,x_7)$ & $(x_5,x_6)$ & $x_7,(x_5,x_7),(x_6,x_7)$\\
\hline
$(x_2,x_5,x_8)$ & $x_5,(x_2,x_5),(x_5,x_8)$ & $(x_2,x_8)$\\
\hline
$(x_3,x_6,x_9)$ & $x_6,(x_3,x_6),(x_6,x_9)$ & $(x_3,x_9)$\\
\hline
$(x_5,x_7,x_8)$ & $(x_5,x_7)$ & $x_8,(x_5,x_8),(x_7,x_8)$\\
\hline
$(x_6,x_7,x_9)$ & $(x_6,x_7)$ & $x_9,(x_6,x_9),(x_7,x_9)$\\
\hline
$(x_7,x_8,x_9)$ & $x_7,(x_7,x_8),(x_7,x_9)$ & $(x_8,x_9)$\\
\hline
$(x_1,x_8,x_9)$ & $(x_8,x_9)$ & $x_1,(x_1,x_8),(x_1,x_9)$\\
\hline
$(x_1,x_2,x_8)$ & $x_8,(x_1,x_8),(x_2,x_8)$ & $(x_1,x_2)$\\
\hline
$(x_1,x_3,x_9)$ & $x_9,(x_1,x_9),(x_3,x_9)$ & $(x_1,x_3)$\\
\hline
\end{tabular}
\\
\caption{The Bongartz (co-)completions with respect to ${\xx}$}
\end{table}

Take $\xx=(x_1,x_2)$ to be a partial cluster and
denote by $\ol{\xx}=(x_1,x_2,x_3)$ and $\mu_{x_3}(\ol{\xx})=(x_1,x_2,x_8)$.
Then by \Cref{def:parBgz}, the partial Bongartz completions with respect to $\xx$ are
\begin{align*}
    \Bgzp{\xx}{x_4}=\Bgzp{\ol{\xx}}{x_4}\cap\Bgzp{\mu_{x_3}(\ol{\xx})}{x_4}=(x_2,x_3,x_4)\cap(x_2,x_4,x_5)=(x_2,x_4),\\
    \Bgzp{\xx}{x_6}=\Bgzp{\ol{\xx}}{x_6}\cap\Bgzp{\mu_{x_3}(\ol{\xx})}{x_6}=(x_3,x_6,x_9)\cap(x_6,x_7,x_9)=(x_6,x_9),\\
    \Bgzp{\xx}{x_9}=\Bgzp{\ol{\xx}}{x_9}\cap\Bgzp{\mu_{x_3}(\ol{\xx})}{x_9}=(x_1,x_3,x_9)\cap(x_1,x_8,x_9)=(x_1,x_9),
\end{align*}
The partial Bongartz (co-)completions could be partial clusters when $\xx$ is not a cluster, which is different from the Bongartz (co)-completion.

Take $\xevo={a x_1+ b x_2}\in\cc(\xx)^\circ$.
By \Cref{def:parallel} and \Cref{def:Xsupp}, we compute their $\xx$-supports as follows
\begin{align*}
    \flowu{\xevo}(x_4)&
    =b x_2 - a x_4&\implies
    \suppu{\xx}(x_4)=(x_2,x_4)\\
    \flowu{\xevo}(x_6)&
    =a x_9 - (a+b) x_6&\implies
    \suppu{\xx}(x_6)=(x_6,x_9)\\
    \flowu{\xevo}(x_9)&
    =(a+b) x_1 - b x_9&\implies
    \suppu{\xx}(x_9)=(x_1,x_9)
\end{align*}
Comparing the partial Bongartz completion and the $\xx$-supports we see the evidence of \Cref{thm:X=xx} in the $A_3$ case.

\begin{figure}
    \centering
    \reviewfigure[width=0.8\linewidth]{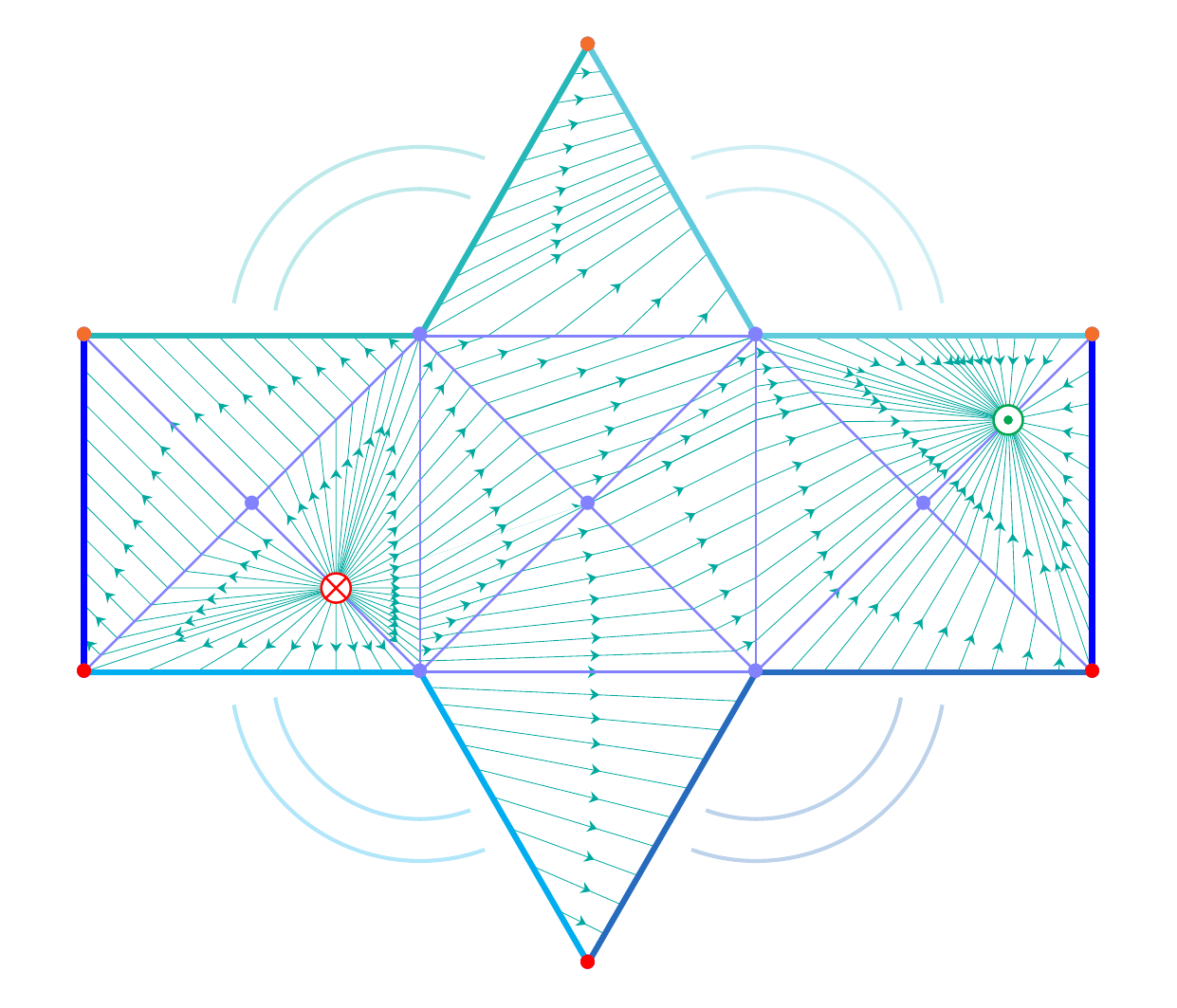}
    \caption{$(\frac{x_1+x_2}{2})$-foliation on $\CCx(\seeds_{A_3})$}
    \label{fig:A3-E1}
\end{figure}

In particular, taking $\xevo=(x_1+x_2)/2$,
the resulting field is illustrated in \Cref{fig:A3-E1}.

Write $\rdom{\xx}(\vv)=\dom{\xx}(\vv)\backslash\xx$ for the reduced dominant sets.
Then the reduced dominant set for each $x_i$ are listed as follows.
\begin{align*}
     &\rdom{\xx}(x_1)=\{x_8,x_9\},\quad
     \rdom{\xx}(x_2)=\{x_4,x_5,x_8\},\quad
     \rdom{\xx}(x_3)=\{x_4,x_6,x_9\},\\
     &\rdom{\xx}(x_4)=\rdom{\xx}(x_5)=\rdom{\xx}(x_6)=\{x_4,x_5,x_6,x_7,x_8,x_9\}=\CCx_1(\seeds)\backslash\xx,\\
     \rdom{\xx}(x_7)&=\{x_5,x_6,x_7,x_8,x_9\},~~
     \rdom{\xx}(x_8)=\{x_5,x_7,x_8,x_9\},~~
     \rdom{\xx}(x_9)=\{x_6,x_7,x_8,x_9\}.
\end{align*}


\end{document}